\documentclass[10pt]{amsart}
\allowdisplaybreaks[4]
\usepackage{amssymb,color}

\usepackage[square,compress,comma, numbers, sort]{natbib}

\usepackage[colorlinks=true, citecolor=blue, linkcolor=blue]{hyperref}
\usepackage{amsthm}

\definecolor{c20}{rgb}{0.,0.7,0.}
\definecolor{c30}{rgb}{0.,0.,1.}
\definecolor{c40}{rgb}{1,0.1,0.7}
\definecolor{c50}{rgb}{1,0,0}
\definecolor{c60}{rgb}{0,0.9,0.1}

\newcommand{\E}[1]{\mathbb{E}\left\{ \left(#1\right)\right\}}

\newcommand{\pk}[1]{\mathbb{P} \left\{#1 \right\} }

\newcommand{\BQN}{\begin{eqnarray}}
\newcommand{\EQN}{\end{eqnarray}}
\newcommand{\BQNY}{\begin{eqnarray*}}
\newcommand{\EQNY}{\end{eqnarray*}}

\newcommand{\BS}{\begin{sat}}
\newcommand{\ES}{\end{sat}}
\newcommand{\BT}{\begin{theo}}
\newcommand{\ET}{\end{theo}}
\newcommand{\BK}{\begin{korr}}
\newcommand{\EK}{\end{korr}}

\newcommand{\BD}{\begin{de}}
\newcommand{\ED}{\end{de}}

\newcommand{\BIT}{\begin{itemize}}
\newcommand{\EIT}{\end{itemize}}
\newcommand{\BDI}{\begin{description}}
\newcommand{\EDI}{\end{description}}

\newcommand{\BRM}{\begin{remarks}}
\newcommand{\ERM}{\end{remarks}}

\newcommand{\BEL}{\begin{lem}}
\newcommand{\EEL}{\end{lem}}

\newtheorem{theo}{Theorem}[section]
\newtheorem{sat}[theo]{Proposition}
\newtheorem{de}[theo]{Definition}
\newtheorem{lem}[theo]{Lemma}

\newtheorem{korr}[theo]{Corollary}

\newtheorem{remarks}[theo]{Remarks}

\newcommand{\la}[1]{\overleftarrow{#1}}

\newcommand{\prooftheo}[1]{ \textbf{Proof of Theorem} \ref{#1} }
\newcommand{\proofprop}[1]{\textbf{Proof of Proposition} \ref{#1}}

\newcommand{\COM}[1]{}

\newcommand{\QED}{\hfill $\Box$ \\}

\def\rw{\rightarrow}

\def\IF{\infty}

\def\rw{\rightarrow}

\def\IF{\infty}

\date{}

\def\rw{\rightarrow}

\def\Var{\text{Var}}

\def\LL{\mathcal{\rho}}
\def\vv{v}

\begin{document}

\title[Gaussian Fields with maximum points of variance on curves ]
{Extremes of Gaussian Random Fields with maximum variance attained over smooth curves}

\author{Peng Liu}

\address{Peng Liu,
 Department of Actuarial Science, University of Lausanne, UNIL-Dorigny 1015 Lausanne, Switzerland
and Mathematical Institute, University of Wroclaw, pl. Grunwaldzki 2/4, 50-384 Wroclaw, Poland}
\email{peng.liu@unil.ch}

\bigskip

\date{\today}
 \maketitle

\bigskip
{\bf Abstract:} Let $X(s,t), (s,t)\in E$, with $E\subset \mathbb{R}^2$ a compact set, be a centered two dimensional Gaussian random field with continuous trajectories and variance function $\sigma(s,t)$. Denote by $\mathcal{L}=\{(s,t): \sigma(s,t)=\max_{(s',t')\in E}\sigma(s',t')\}$. In this contribution, we derive the exact asymptotics of $\pk{\sup_{(s,t)\in E}X(s,t)>u}$, as $u\rw\IF$, under condition that $\mathcal{L}$ is a smooth curve.  We illustrate our findings by an application concerned with extremes of the aggregation of two independent fractional Brownian motions.
 \\



\section{Introduction}

Tail asymptotics of supremum of Gaussian processes and Gaussian random fields are investigated substantially in the literature, most of which consider the stationary Gaussian random  fields or non-stationary case with unique maximum point of variance,  see e.g., \cite{PicandsA, Pit72, berman1987extreme, Pit96, Michna2, Michna1, DE2002, DI2005, DE2014, DEJ14, Pit20,Tabis,Marek,MR2462286,MR3583761,MR0494458}. Recently, in  \cite{KEP20162}, the extremes of Gaussian random fields with unique point of maximum variance and  more general local variance and correlation structure have been considered.\\
 Specifically, let $X(s,t), (s,t)\in E=[-T_1,T_1]\times [-T_2,T_2]$ be a Gaussian random field with continuous trajectories and variance $\sigma(s,t)$ attaining the maximum at unique point $(0,0)$ with $\sigma(0,0)=1$ and set
$$A=\left(\begin{array}{cc}
a_{11} & a_{12}\\
a_{21}& a_{22}\\
\end{array}\right), \quad B=\left(\begin{array}{cc}
b_{11}&b_{12}\\
b_{21}&b_{22}\\
\end{array}\right).$$
 \cite{KEP20162} investigates the Gaussian random field $X$ with regularly varying dependence structure, i.e., for $A\neq 0$, $|B|\neq 0$,
\BQN\label{cor}
1-r(s,t,s_1,t_1)\sim \rho_1^2(|a_{11}(s-s_1)+a_{12}(t-t_1)|)+\rho_2^2(|a_{21}(s-s_1)+a_{22}(t-t_1)|),
\EQN
holds for $|s-s_1|, |t-t_1|, |B(s,t)^\top|\rw 0$,
and
\BQN\label{var}
1-\sigma(s,t)\sim v_1^2(|b_{11}s+b_{12}t|)+v_2^2(|b_{21}s+b_{22}t|), \quad |B(s,t)^\top|\rw 0,
\EQN
where $|B|$ is the determinant of $B$, $D^\top$ is the transpose of $D$ for any matrix $D$, $\rho_i, v_i\geq 0$, $\rho_i\in \mathcal{R}_{\alpha_i/2},\alpha_i\in (0,2]$ and  $v_i\in \mathcal{R}_{\beta_i/2}, \beta_i>0$ (here $\mathcal{R}_\gamma$ stands for the class of regularly varying functions at 0 with index $\gamma$, see, e.g., \cite{BI1989}). Due to the regularity of $\rho_i, v_i, i=1,2,$ and the locally non-additive dependence of correlation and variance structure, qualitatively new types of tail asymptotics for $M=\sup_{(s,t) \in E} X(s,t)$ have been obtained in \cite{KEP20162}.
It is worthwhile to note that the invertibility of $B$ plays the key role to ensure that the maximum of $\sigma$ is attained at a unique point. Correspondingly, if $B$ is non-invertible, (\ref{var}) shows that $\sigma$ attains maximum at $\{(s,t): B(s,t)^\top=0\}$, which is an uncountable set. The non-uniqueness of the maximum variance renders the the problem of the tail asymptotics of $M$ essentially quite different from the case considered in \cite{KEP20162}. This has already been addressed in \cite{nonhomoANN} for standard assumptions on the variance and the covariance of the Gaussian random field $X$, i.e., $v_i,r_i$ are simple power functions. In this paper, as a continuation and complement to the results obtained in  \cite{KEP20162,nonhomoANN}, we investigate the case of $|B|=0$ and $B\neq 0$ considering general structures for variance and covariance function of $X$.

\COM{Some partial results on this problem were considered in \cite{nonhomoANN} with the Gaussian fields having the following local structure:
$$1-r(s,t,s',t')\sim |a_{11}(s-s_1)|^{\alpha_1}+|a_{21}(s-s_1)+a_{22}(t-t_1|)|^{\alpha_2},  \quad |s-s_1|, |t-t_1|\rw 0, $$
and
$$1-\sigma(s,t)\sim b_1|t-T|^\beta, |t-T|\rw 0,$$
where $a_{11}>0, a_{21}\geq 0, a_{22}\neq 0, b_1>0$.\\
}
As motivated in \cite{nonhomoANN}, Gaussian random fields with non-unique point of maximum (attained on a line) appear in connection with the  Shepp statistics defined by
$$Y(t)=\sup_{s\in [0,T_1]}X(s+t)-X(s), \quad t\in [0,T_2],$$
where $\{X(t), t\geq 0\}$ is a centered Gaussian process with continuous trajectories.
Of interest in statistics is the tail asymptotics  of the supremum of $Y$, i.e., the asymptotics of
\BQN\label{Shepp}\pk{\sup_{t\in [0,T_2]}Y(t)>u}=\pk{\sup_{(s,t)\in [0,T_1]\times[0,T_2]}(X(s+t)-X(s))>u}, \ \ u\rw\IF,
\EQN
where the Gaussian random field $X(s+t)-X(s), (s,t)\in [0,T_1]\times[0,T_2]$ satisfies (\ref{cor}) and (\ref{var}) with $B\neq 0$ and $|B|=0$ and the maximum points of variance  are $\{(s,T): 0\leq s\leq T_1\}$. In \cite{nonhomoANN} the asymptotics of \eqref{Shepp} is derived for $X$ being an fBm, Ornstein-Uhlenbeck process and Brownian bridge. We claim that the asymptotics of \eqref{Shepp} for general Gaussian process $X$ can be derived by using our results in this paper. Since our main focus is on the curve case,  in this paper we don't present the results. \\
As mentioned above, \cite{nonhomoANN} investigates Gaussian random fields where the maximum point of the variance is attained on a line. In this paper we shall consider more general cases where the maximum point of the variance of $X$ is attained on a smooth curve. We give an application of this new result to the study of the aggregation of two independent fractional Brownian motions.  \\
Brief outline of the rest of the paper: Section 2 is dedicated to the case that the  variance of $X$ is maximal on a line followed by Section 3 which extends those findings to the case that the maximum of the variance of $X$ is attained on a smooth curve. An application is displayed in Section 4 followed by Section 5 which contains all the proofs.

\section{Maximum variance attained over a line }
Let $X(s,t), (s,t)\in E=[-T_1,T_1]\times[-T_2,T_2]$ be a centered Gaussian field with continuous trajectories and correlation and variance functions satisfying (\ref{cor}) and (\ref{var}) with $A, B\neq 0, |B|=0$.  We study the asymptotics of
\BQNY
\pi(u)=\pk{\sup_{(s,t)\in E }X(s,t)>u}, \ \ u\rw\IF.
\EQNY
Let in the following $B_{\alpha}(t), t\geq 0$ be a standard fractional Brownian motion (fBm) with Hurst index $\alpha\in (0,2]$  and covariance satisfying
$$Cov\left(B_{\alpha}(s), B_{\alpha}(t)\right)=\frac{|s|^{\alpha}+|t|^{\alpha}-|t-s|^{\alpha}}{2},$$
and define
$$W_{\alpha_1,\alpha_2}(s,t)=\sqrt{2}B_{\alpha_1}(s)+\sqrt{2}B_{\alpha_2}(t)-|s|^{\alpha_1}-|t|^{\alpha_2},$$ with
$B_{\alpha_i}, i=1,2$ two independent fBms with indices $\alpha_1,\alpha_2\in (0,2]$.\\
In order to proceed with our findings, we give first the definition of Pickands and Piterbarg constants which are crucially important in the extreme theory of Gaussian processes and fields.
Define  first the Pickands constant by
$$\mathcal{H}_{\alpha}[0,S]=\E{\sup_{t\in [0,S]}e^{\sqrt{2}B_{\alpha}(t)-|t|^{\alpha}}}, \quad \mathcal{H}_{\alpha}=\lim_{S\rw\IF}\frac{ \mathcal{H}[0,S]}{S},$$
and Piterbarg constant by
$$\mathcal{P}_{\alpha}^\gamma[S_1,S_2]=\E{\sup_{t\in [S_1,S_2]}e^{\sqrt{2}B_{\alpha}(t)-(1+\gamma)|t|^{\alpha}}}, \quad  \mathcal{P}_{\alpha}^\gamma=\lim_{S\rw\IF}\mathcal{P}_{\alpha}^\gamma[S,S], \quad  \widehat{\mathcal{P}}_{\alpha}^\gamma=\lim_{S\rw\IF} \mathcal{P}_\alpha^\gamma[0,S],$$
for $\alpha\in (0,2]$, $\gamma>0, S_1<S_2$.  We denote by $\mathcal{P}_{\alpha}^\IF=\widehat{\mathcal{P}}_{\alpha}^\IF=1$.
Moreover, let
$$\mathcal{H}_{\alpha_1,\alpha_2}^{\gamma, b,\beta}(S)=\mathbb{E}\left(\sup_{\{|s+bt|\leq S, 0\leq t\leq S\}}e^{W_{\alpha_1,\alpha_2}(s,t)-\gamma|s+bt|^\beta}\right), \quad \widehat{\mathcal{H}}_{\alpha_1,\alpha_2}^{\gamma, b,\beta}(S)=\mathbb{E}\left(\sup_{\{0\leq s+bt\leq S, 0\leq t\leq S\}}e^{W_{\alpha_1,\alpha_2}(s,t)-\gamma|s+bt|^\beta}\right),$$
where  $\beta>0, b\in\mathbb{R}, S\geq 0.$ For simplicity, denote by
\BQN\label{piterbarg}
\mathcal{H}_{\alpha}^{\gamma,b}(S):=\mathcal{H}_{\alpha,\alpha}^{\gamma, b,\alpha}(S), \quad\mathcal{H}_{\alpha}^{\gamma,b}:=\lim_{S\rw\IF}\frac{\mathcal{H}_{\alpha}^{\gamma,b}(S)}{S}, \quad \widehat{\mathcal{H}}_{\alpha}^{\gamma,b}(S):=\widehat{\mathcal{H}}_{\alpha,\alpha}^{\gamma, b,\alpha}(S), \quad\widehat{\mathcal{H}}_{\alpha}^{\gamma,b}:=\lim_{S\rw\IF}\frac{\widehat{\mathcal{H}}_{\alpha}^{\gamma,b}(S)}{S}.
\EQN
For the extensions and related properties of Pickands-Piterbarg constants, see e.g., \cite{PicandsA, Pit72, Pit96, Michna2,Michna1, DE2002, DI2005, nonhomoANN, Krzys2006Pickands, PitKonst,
	debicki2008note,  DiekerY,
Harper1, Harper2,DE2014, Pit20, DM, SBK} and the references therein. One can refer to \cite{nonhomoANN} and \cite{KEP20162} for  the existence of the limit in (\ref{piterbarg}). \\
 Throughout this paper we shall assume that rank$(A)=2$, the case rank$(A)=1$ shall be considered separately due to too many technical details. \\
 Let $Z(s,t)=X(A^{-1}(s,t)^\top)$. Then by (\ref{cor}) and (\ref{var1}), we have
\BQN\label{cor1}
1-r_Z(s,t,s_1,t_1)\sim \rho_1^2(|s-s_1|)+\rho_2^2(|t-t_1|), \quad |s-s_1|, |t-t_1|, |BA^{-1}(s,t)^\top|, |BA^{-1}(s_1,t_1)^\top|\rw 0,
\EQN
and
\BQN\label{var1}
1-\sigma_Z(s,t)\sim v_1^2(|c_{11}s+c_{12}t|)+v_2^2(|c_{21}s+c_{22}t|),
\EQN
with
$$BA^{-1}=\left(\begin{array}{cc}c_{11} & c_{12}\\
c_{21}& c_{22}\end{array}\right).$$
Note that $BA^{-1}\neq 0$ and $|BA^{-1}|=0$. Without loss of generality, we assume that $c_{11}\neq 0$. Then it follows that (\ref{cor1}) and (\ref{var1}) can be rewritten as
\BQN\label{cor2}
1-r(s,t,s_1,t_1)\sim \rho_1^2(|s-s_1|)+\rho_2^2(|t-t_1|), \quad |s-s_1|, |t-t_1|, |s+bt|, |s_1+bt_1|\rw 0,
\EQN
and
\BQN\label{var2}
1-\sigma(s,t)\sim v^2(|s+bt|), \quad |s+bt|\rw 0,
\EQN
with $b=\frac{c_{12}}{c_{11}}$ and
$$v=\sqrt{|c_{11}|^{\beta_1}v_1^2+|c_{21}|^{\beta_2}v_2^2}\in \mathcal{R}_{\beta/2}, \quad \beta=\min(\beta_1, \beta_2)I_{\{c_{21}\neq 0\}}+\beta_1I_{\{c_{21}=0\}}.$$
 For further analysis, we assume that
\BQN\label{ratio}
\lim_{t\downarrow 0}\frac{\rho_2^2(t)}{\rho_1^2(t)}=\eta\in [0,\IF], \quad \lim_{t\downarrow 0}\frac{v^2(t)}{\rho_i^2(t)}=\gamma_i\in [0,\IF].
\EQN
 Let $X(s,t), (s,t)\in E$ be a centered Gaussian random field with continuous trajectories. Suppose $\sigma(s,t)$ attains its maximum , which equals $1$, at $\mathcal{L}:=\{(s,t)\in [-T_1, T_1]\times[-T_2,T_2], |s+bt|=0\}$. Moreover, assume that
 \BQN\label{local}
 Corr\left(X(s,t), X(s',t')\right)<1, \quad (s,t)\neq (s',t'), (s,t), (s',t')\in [-T_1,T_1]\times [-T_2,T_2].
 \EQN
Hereafter, let $\Psi(u)$ denote the tail distribution of a standard normal random variable. For any $v\in \mathcal{R}_\lambda, \lambda>0$, $\overleftarrow{v}$ denotes the asymptotic (unique) inverse of $f\in \mathcal{R}_\gamma$ defined by
$\overleftarrow{v}(x)=\inf\{y\in (0,1]: v(y)>x\}, \quad x>0.$
See, e.g., \cite{BI1989} for the definitions and properties of asymptotic inverse functions. Recall that we denote by $\mathcal{P}_{\alpha}^\IF=1$.\\
 $\diamond$  {\underline{Case 1.}} $b=0$.
 \BT\label{Th1} Suppose that (\ref{cor2})-(\ref{local})  hold with $b=0$.\\
 If $\gamma_1=0$, then
$$\pi(u)\sim 4T_2 \Gamma(1/\beta+1)\prod_{i=1}^2\mathcal{H}_{\alpha_i}\frac{\overleftarrow{\vv}(1/u)}{\overleftarrow{\LL}_1(1/u)\overleftarrow{\LL}_2(1/u)}\Psi(u).$$
If $\gamma_1\in (0,\IF]$, then
$$\pi(u)\sim  2T_2\mathcal{H}_{\alpha_2}\mathcal{P}_{\alpha_1}^{\gamma_1}\frac{1}{\overleftarrow{\LL}_2(1/u)}\Psi(u).$$
\ET
$\diamond$  {\underline{Case 2.}} $b\neq 0$ and $\eta\in (0,\IF)$.
\BT\label{Th2} Suppose that (\ref{cor2})-(\ref{local})  hold with $b\neq 0$ and further assume that $\eta\in (0,\IF)$.\\
i) If $\gamma_1=0$, then
$$\pi(u)\sim 4\min(T_2, T_1/|b|)\Gamma(1/\beta+1)\prod_{i=1}^2\mathcal{H}_{\alpha_i}\frac{\overleftarrow{\vv}(1/u)}{\overleftarrow{\LL}_1(1/u)\overleftarrow{\LL}_2(1/u)}\Psi(u).$$
If $\gamma_1\in (0,\IF)$, then
$$\pi(u)\sim \mathcal{H}_{\alpha_1}^{\gamma_1, b\eta^{-1/\alpha_1}}\frac{2\min(T_2, T_1/|b|)}{\la{\rho}_2(u^{-1})}\Psi(u).$$
If $\gamma_1=\IF$, then
$$\pi(u)\sim (|b|^{\alpha_1}\eta^{-1}+1)^{1/\alpha_1}\mathcal{H}_{\alpha_1}\frac{2\min(T_2, T_1/|b|)}{\la{\rho}_2(u^{-1})}\Psi(u).$$
\ET
\BRM\label{remark0} Assume that $X(s,t), (s,t)\in E$ is a Gaussian random field with $E=([S_1,S_2]\times [T_1,T_2])\cap \{(s,t): s_0+s+b(t+t_0)\geq 0\}$ and $(s_0, t_0)\in (S_1,S_2)\times (T_1,T_2)$. The maximum of variance is attained over
$$E\cap \{(s,t): s_0+s+b(t+t_0)=0\}=\{(s,t): s_0+s+b(t+t_0)=0, t_1\leq t\leq t_2\}.$$ This means that the maximum  is attained at the boundary of $E$. If (\ref{cor2})-(\ref{local}) are satisfied with $[-S_1,S_1]\times [-T_1, T_1]$ replaced by $E$ and $s+bt$ replaced by $s_0+s+b(t+t_0)$, then Theorem \ref{Th1} holds with $ 2T_2$ replaced by $t_2-t_1$, $\Gamma(1/\beta+1)$ replaced by $\frac{1}{2}\Gamma(1/\beta+1)$ and  $\mathcal{P}_{\alpha_1}^{\gamma_1}$ replaced by $\widehat{\mathcal{P}}_{\alpha_1}^{\gamma_1}$. Theorem \ref{Th2} holds with  $ 2\min(T_2, T_1/|b|)$ replaced by $t_2-t_1$, $\Gamma(1/\beta+1)$ replaced by $\frac{1}{2}\Gamma(1/\beta+1)$ and $\mathcal{H}_{\alpha_1}^{\gamma_1, b\eta^{-1/\alpha_1}}$ replaced by $\widehat{\mathcal{H}}_{\alpha_1}^{\gamma_1, b\eta^{-1/\alpha_1}}$.
\ERM
\subsubsection{Discussion} In this subsection, we will show that all the cases for rank$(A)=2$ can be reduced to Case 1-2. Without loss of generality, we assume that $c_{11}\neq 0$. Then we continue with the analysis from (\ref{cor2}) and (\ref{var2}).
The following hold:
\begin{itemize}
\item[$\diamond$] $\underline{b=0}$: this case is covered by Case 1.
\item[$\diamond$] $\underline{b\neq 0,\eta=0}$: let $Y(s,t)=X(t, \frac{s-t}{b}), (s,t)\in E_Y=\{(s,t), (t, \frac{s-t}{b})\in [-T_1,T_1]\times [-T_2,T_2]\}$. Then $\sigma_Y$ attains its maximum on $\{(0, t), |t|\leq \min(|T_1|,|bT_2|)\}$ and in light of Lemma 6.4 in \cite{KEP20162}, we have
$$
    1-r_Y(s,t,s_1,t_1)\sim \rho_1^2(|t-t_1|)+\rho_2^2(|b^{-1}(s-s_1-(t-t_1))|)\sim |b|^{-\alpha_2}\rho_2^2(|s-s_1|)+\rho_1^2(|t-t_1|),
  $$
for $\quad |s-s_1|, |t-t_1|, |s|\rw 0,$ and
$1-\sigma_Y(s,t)\sim v^2(|s|),  s\rw 0$.
\item[$\diamond$] $\underline{b\neq 0, \eta\in (0,\IF)}$: this case is covered by Case 2.
\item[$\diamond$] $\underline{b\neq 0, \eta=\IF}$:
 let $Y(s,t)=X(s-bt, t), (s,t)\in E_Y=\{(s,t), (s-bt, t)\in [-T_1,T_1]\times [-T_2,T_2]\}$. Then $\sigma_Y$ attains its maximum over $\{(0, t), |t|\leq \min\left(|T_1/b|, |T_2|\right)\}$ and in light of Lemma 6.4 in \cite{KEP20162}, we have
$$
    1-r_Y(s,t,s_1,t_1)\sim \rho_1^2(|s-s_1-b(t-t_1)|)+\rho_2^2(|t-t_1|)\sim \rho_1^2(|s-s_1|)+\rho_2^2(|t-t_1|),
  $$
for $\quad |s-s_1|, |t-t_1|, |s|\rw 0,$ and
$1-\sigma_Y(s,t)\sim v^2(|s|), \quad s\rw 0$.
\end{itemize}

\section{ Maximum  variance attained over a curve }
In the previous section, the maximum of variance is attained over a straight line. Whereas, in this section we focus on the scenarios that the maximum of variance is attained over a smooth curve. Let $X(s,t), (s,t)\in E=[S_1,S_2]\times[T_1,T_2]$ be a centered Gaussian process with continuous trajectories, maximum variance equaling $1$ attained over $\mathcal{L}=\{(s,t): (s,t)\in E, s=f(t)\}$ with $f$ a function. Let $f, g$ be two continuous functions with $g$ satisfying $0<c_1\leq g(t)\leq c_2<\IF, t\in [T_1, T_2]$ and $f$ satisfying \\

{\bf F}. $f\in C^1((T_1,T_2))$ and  $f(T_i)\in [S_1,S_2], i=1,2$. Moreover, $\inf_{t\in (T_1,T_2)}|f'(t)|>0$.\\

Denote by $\la{f}(t)$ the inverse function of $f$. Inspired by (\ref{cor2}) and (\ref{var2}) and the application in section 4, we assume that
\BQN\label{cor4}
\lim_{\delta\rw 0}\sup_{(s,t), (s',t')\in E, (s,t)\neq (s',t'),|s-s'|, |t-t'|, |s-f(t)|, |s'-f(t')|<\delta}\left|\frac{1-r(s,t,s',t')}{\rho_1^2(|s-s'|)+\rho_2^2(|t-t'|)}-1\right|=0,
\EQN
and
\BQN\label{var4}
\lim_{\delta\rw 0}\sup_{(s,t)\in E, |s-f(t)|<\delta}\left|\frac{1-\sigma(s,t)}{ v^2(g(t)|s-f(t)|)}-1\right|=0.
\EQN
Moreover, assume that
 \BQN\label{local1}
 Corr\left(X(s,t), X(s',t')\right)<1, \quad (s,t)\neq (s',t'), (s,t), (s',t')\in [S_1,S_2]\times [T_1,T_2].
 \EQN
Note that    assumption  (\ref{var4}) can be justified by (\ref{var5}) and (\ref{var6}), which are the local behaviors of variances of the crucial example in section 4.  We use the same natation as in (\ref{ratio}) in the following theorems.
\BT\label{th31}
Suppose that (\ref{cor4})-(\ref{local1}) hold with $f$ satisfying {\bf F} and $\eta =0$.\\
If $\gamma_2=0$, then
$$\pi(u)\sim 2\int_{T_1}^{T_2}(g(t))^{-1}dt\Gamma(1/\beta+1)\prod_{i=1}^2\mathcal{H}_{\alpha_i}\frac{\overleftarrow{\vv}(1/u)}{\overleftarrow{\LL}_1(1/u)
\overleftarrow{\LL}_2(1/u)}\Psi(u).$$
If $\gamma_2\in (0,\IF]$, then
$$\pi(u)\sim \mathcal{H}_{\alpha_1}\int_{T_1}^{T_2} \mathcal{P}_{\alpha_2}^{\gamma_2|f'(t)g(t)|^{\alpha_2}}|f'(t)|dt\frac{\Psi(u)}{\la{\rho}_1(1/u)}.$$
\ET
\BT\label{th32}
Suppose that (\ref{cor4})-(\ref{local1}) hold with $f$ satisfying {\bf F} and $\eta \in (0,\IF)$.\\
If $\gamma_1=0$, then
$$\pi(u)\sim 2\int_{T_1}^{T_2}(g(t))^{-1}dt\Gamma(1/\beta+1)\prod_{i=1}^2\mathcal{H}_{\alpha_i}\frac{\overleftarrow{\vv}(1/u)}{\overleftarrow{\LL}_1(1/u)\overleftarrow{\LL}_2(1/u)}\Psi(u).$$
If $\gamma_1\in (0,\IF)$, then
$$\pi(u)\sim \int_{T_1}^{T_2} \mathcal{H}_{\alpha_1}^{\gamma_1(g(t))^{\alpha_1}, -\eta^{-1/\alpha_1}|f'(t)|}dt\frac{\Psi(u)}{\la{\rho}_2(1/u)}.$$
\ET
If $\gamma_1=\IF$, then
$$\pi(u)\sim \int_{T_1}^{T_2}(\eta^{-1}|f'(t)|^{\alpha_1}+1)^{1/\alpha_1}dt\mathcal{H}_{\alpha_1}\frac{\Psi(u)}{\la{\rho}_2(1/u)}.$$
\BT\label{th33}
Suppose that (\ref{cor4})-(\ref{local1}) hold with $f$ satisfying {\bf F} and $\eta =\IF$.\\
If $\gamma_1=0$, then
$$\pi(u)\sim 2\int_{T_1}^{T_2}(g(t))^{-1}dt\Gamma(1/\beta+1)\prod_{i=1}^2\mathcal{H}_{\alpha_i}\frac{\overleftarrow{\vv}(1/u)}{\overleftarrow{\LL}_1(1/u)\overleftarrow{\LL}_2(1/u)}\Psi(u).$$
If $\gamma_1\in (0,\IF]$, then
$$\pi(u)\sim \mathcal{H}_{\alpha_2}\int_{T_1}^{T_2}\mathcal{P}_{\alpha_1}^{\gamma_1(g(t))^\beta}dt\frac{\Psi(u)}
 {\la{\rho}_2(1/u)}.$$
\ET
\BRM\label{remark} Assume that $X(s,t), (s,t)\in E$ is a Gaussian field with $E^{\pm}=([S_1,S_1]\times [T_1,T_2])\cap \{(s,t): \pm(s-f(t))\geq 0\}$. The maximum of variance is attained over the curve $([S_1,S_1]\times [T_1,T_2])\cap \{(s,t): s-f(t)=0\}$. This means that the maximum of variance is attained at the boundary of $E^\pm$. We further assume that (\ref{cor4})-(\ref{local1}) and {\bf F} are all satisfied. Then Theorems \ref{th31}-\ref{th33} hold with $\Gamma(1/\beta+1)$ replaced by $\frac{1}{2} \Gamma(1/\beta+1)$, $\mathcal{P}_{\alpha}^{\gamma}$ by $\widehat{\mathcal{P}}_{\alpha}^{\gamma}$ and $\mathcal{H}_{\alpha}^{\gamma, b}$ by $\widehat{\mathcal{H}}_{\alpha}^{\gamma, b}$.
\ERM
\section{Applications}

Let $B_{\alpha_i}(t), i=1,2$ be independent fBms with indices $\alpha_i\in (0,2), i=1,2$ respectively. Of interest is the  asymptotic of
$$\pk{\sup_{(s,t)\in E_{\alpha_1,\alpha_2}}(B_{\alpha_1}(s)+B_{\alpha_2}(t))>u}, \quad \hbox{with }E_{\alpha_1,\alpha_2}=\{(s,t): |s|^{\alpha_1}+|t|^{\alpha_2}\leq 1\}.$$
Note that $\sigma(s,t)=\sqrt{\Var(B_{\alpha_1}(s)+B_{\alpha_2}(t))}$ attains $1$ at
 $\mathcal{L}=\{(s,t): |s|^{\alpha_1}+|t|^{\alpha_2}=1\}$.
Applying Theorems \ref{th31}-\ref{th33}, we derive the following results.
\BS\label{cor11} Assume that $\alpha_2>\alpha_1$. If $\alpha_2<1$, then
\BQNY
\pk{\sup_{(s,t)\in E_{\alpha_1,\alpha_2}}(B_{\alpha_1}(s)+B_{\alpha_2}(t))>u}\sim \frac{2^{3-1/\alpha_1-1/\alpha_2}}{\alpha_1}\prod_{i=1}^2\mathcal{H}_{\alpha_i}\int_{0}^{1}(1-t^{\alpha_2})^{1/\alpha_1-1}dt
u^{2/\alpha_1+2/\alpha_2-2}\Psi(u).
\EQNY
If $\alpha_2=1$, then
\BQNY
\pk{\sup_{(s,t)\in E_{\alpha_1,\alpha_2}}(B_{\alpha_1}(s)+B_{\alpha_2}(t))>u}\sim 2^{3-1/\alpha_1}\mathcal{H}_{\alpha_1}u^{2/\alpha_1}\Psi(u).
\EQNY
If $\alpha_2>1$, then
\BQNY
\pk{\sup_{(s,t)\in E_{\alpha_1,\alpha_2}}(B_{\alpha_1}(s)+B_{\alpha_2}(t))>u}\sim \frac{ 2^{2-1/\alpha_1}\alpha_2}{\alpha_1}\mathcal{H}_{\alpha_1}\int_{0}^{1}(1-t^{\alpha_2})^{1/\alpha_1-1}t^{\alpha_2-1}dt u^{2/\alpha_1}\Psi(u).
\EQNY
\ES
\BS\label{cor12} Assume that $\alpha_1=\alpha_2=\alpha$. If $\alpha<1$, then
\BQNY
\pk{\sup_{(s,t)\in E_{\alpha_1,\alpha_2}}(B_{\alpha_1}(s)+B_{\alpha_2}(t))>u}\sim \frac{2^{3-2/\alpha}}{\alpha}\left(\mathcal{H}_{\alpha}\right)^2\int_{0}^{1}(1-t^{\alpha})^{1/\alpha-1}dt
u^{4/\alpha-2}\Psi(u).
\EQNY
If $\alpha=1$, then
\BQNY
\pk{\sup_{(s,t)\in E_{\alpha_1,\alpha_2}}(B_{\alpha_1}(s)+B_{\alpha_2}(t))>u}\sim 2\widehat{\mathcal{H}}_{1}^{1,-1}
u^{2}\Psi(u).
\EQNY
If $\alpha>1$, then
\BQNY
\pk{\sup_{(s,t)\in E_{\alpha_1,\alpha_2}}(B_{\alpha_1}(s)+B_{\alpha_2}(t))>u}\sim 2^{2-1/\alpha}\mathcal{H}_{\alpha}\int_0^1\left(1+t^{\alpha(\alpha-1)}(1-t^{\alpha})^{1-\alpha}\right)^{1/\alpha}dt u^{2/\alpha}\Psi(u).
\EQNY
\ES

\section{Proofs}
Throughout this section, we denote by $\mathbb{Q}$ a positive constant that may be different from line to line.\\
 \prooftheo{Th2} For the analysis of the asymptotics, we introduce the following notation:
\BQNY
E(u)&=&\{(s,t)\in E, |s+bt|\leq \la{v}(\ln u/u)\}, \quad \quad I_{k}(u)=[k\la{\rho}_1(u^{-1})S, (k+1)\la{\rho}_1(u^{-1})S],\\
J_l(u)&=& [l\la{\rho}_2(u^{-1})S, (l+1)\la{\rho}_2(u^{-1})S], \quad I_{k,l}(u)=I_{k}(u)\times J_l(u), \\
N_1(u)&=&\left[\frac{\la{v}(\ln u/u)}{\la{\rho}_1(u^{-1})S}\right], \quad N_2(u)=\left[\frac{\min(T_2, T_1/|b|)}{\la{\rho}_2(u^{-1})S}\right],\\
u_{k,l}^{+\epsilon}&=&u\left(1+(1+\epsilon)\sup_{(s,t)\in I_{k,l}(u)}v^2(|s+bt|)\right), \quad u_{k,l}^{-\epsilon}=u\left(1+(1-\epsilon)\inf_{(s,t)\in I_{k,l}(u)}v^2(|s+bt|)\right),\\
 K_u^+&=&\{(k,l):I_{k,l}(u)\cap E(u)\neq \emptyset\}, \quad K_u^-=\{(k,l):I_{k,l}(u)\subset E(u)\}, \\
 E_l^+&=&\{k: I_{k,l}(u)\cap E(u)\neq \emptyset \}, \quad E_l^-=\{k: I_{k,l}(u)\subset E(u) \}, \\
  \Lambda_1(u)&=&\{(k,l,k_1,l_1): (k,l), (k_1,l_1)\in K_u^-, (k,l)\neq (k_1,l_1), k\leq k_1, I_{k,l}\cap I_{k_1,l_1}(u)\neq \emptyset\},\\
\Lambda_2(u)&=&\{(k,l,k_1,l_1): (k,l), (k_1,l_1)\in K_u^-, , k\leq k_1, I_{k,l}\cap I_{k_1,l_1}(u)= \emptyset, |\la{\rho}_2(u^{-1})(l-l_1)S|\leq \epsilon \},\\
 \Lambda_3(u)&=&\{(k,l,k_1,l_1): (k,l), (k_1,l_1)\in K_u^-, , k\leq k_1, I_{k,l}\cap I_{k_1,l_1}(u)= \emptyset, |\la{\rho}_2(u^{-1})(l-l_1)S|\geq \epsilon  \},
 \EQNY
 with $\epsilon>0$ sufficiently small.\\
 It follows that for $\epsilon_1$ sufficiently small,
 \BQN\label{main2}
 \pi_1(u)\leq\pi(u)\leq \pi_1(u)+\pk{\sup_{E\setminus E_{\epsilon_1}}X(s,t)>u}+\pk{\sup_{E_{\epsilon_1}\setminus E(u)}X(s,t)>u},
 \EQN
 with
 $$\pi_1(u)=\pk{\sup_{(s,t)\in E(u)}X(s,t)>u}, \quad  E_{\epsilon_1}=\{(s,t):|s+bt|\leq \epsilon_1\}\cap E.$$
By the fact that $$\sup_{(s,t)\in E\setminus E_{\epsilon_1}}\sigma(s,t)<1-\delta$$ with $0<\delta<1$ and using Borell-TIS inequality (\cite{AdlerTaylor}), we have
\BQN\label{eq01}
\pk{\sup_{(s,t)\in E\setminus E_{\epsilon_1}}X(s,t)>u}\leq e^{-\frac{(u-a)^2}{2(1-\delta)^2}}
\EQN
with $a=\mathbb{E}\left(\sup_{(s,t)\in E\setminus E_{\epsilon_1}}X(s,t)\right)$. In light of (\ref{cor2}) and (\ref{var2}), we have for $u$ sufficiently large,
\BQNY
\sup_{(s,t)\in  E\setminus E(u)}\sigma(s,t)\leq 1-\left(\frac{\ln u}{2u}\right)^2,
\EQNY
and, with $\epsilon_2$ sufficiently small, for $|s-s_1|, |t-t_1|, |s+bt|, |s_1+b t_1|\leq \epsilon_2$
\BQNY
\mathbb{E}\left((\overline{X}(s,t)-\overline{X}(s_1,t_1))^2\right)=2(1-r(s,t,s_1,t_1))&\leq& 4(\rho_1^2(|s-s_1|)+\rho_2^2(|t-t_1|))\\
&\leq& \mathbb{Q}\left(|s-s_1|^{\alpha_1/2}+|t-t_1|^{\alpha_2/2}\right),
\EQNY
implying that
$$\mathbb{E}\left((\overline{X}(s,t)-\overline{X}(s_1,t_1))^2\right)\leq \mathbb{Q} \left(|s-s_1|^{\alpha_1/2}+|t-t_1|^{\alpha_2/2}\right), (s,t), (s_1,t_1)\in E_{\epsilon_1}\setminus E(u).$$
Consequently, by Lemma 5.1 in \cite{KEP2016} or Theorem 8.1 in \cite{Pit96}, for $u$ large enough,
\BQN\label{eq2}
\pk{\sup_{(s,t)\in E_{\epsilon_1}\setminus E(u)}X(s,t)>u}&\leq& \pk{\sup_{(s,t)\in E_{\epsilon_1}\setminus E(u)}\overline{X}(s,t)>\frac{u}{\sqrt{1-\left(\frac{\ln u}{2u}\right)^2}}}\nonumber\\
&\leq& \mathbb{Q}u^{4/\alpha_1+4/\alpha_2}\Psi\left(\frac{u}{1-\left(\frac{\ln u}{2u}\right)^2}\right),
\EQN
which together with (\ref{main2}), (\ref{eq01}) and the fact that $\pi_1(u)\geq \Psi(u)$ leads to
\BQNY
\pi(u)\sim \pi_1(u), \ \ u\rw\IF.
\EQNY
 Next we focus on $\pi_1(u)$. Without loss of generality, assume that $|b|T_2\leq T_1$.\\
 \underline{$\gamma_1=0$}. Bonferroni inequality leads to
 \BQN\label{main1}
 \sum_{(k,l)\in K_u^-}\pk{\sup_{(s,t)\in I_{k,l}(u)}\overline{X}(s,t)>u_{k,l}^{+\epsilon}}-\sum_{i=1}^3\Sigma_i(u)\leq \pi_1(u)\leq \sum_{(k,l)\in K_u^+}\pk{\sup_{(s,t)\in I_{k,l}(u)}\overline{X}(s,t)>u_{k,l}^{-\epsilon}},
 \EQN
 where
 $$\Sigma_i(u)=\sum_{(k,l,k_1,l_1)\in \Lambda_i(u)}\pk{\sup_{(s,t)\in I_{k,l}(u)}\overline{X}(s,t)>u_{k,l}^{-\epsilon}, \sup_{(s,t)\in I_{k_1,l_1}(u)}\overline{X}(s,t)>u_{k_1,l_1}^{-\epsilon}}, \quad i=1,2,3.$$
 In order to apply Lemma \ref{PIPI}, we set
 $$\quad X_{u,k,l}(s,t)=\overline{X}(k\la{\rho}_1(u^{-1})S+s, l\la{\rho}_2(u^{-1})S+t), (s,t)\in I_{0,0}(u), (k,l)\in K_u^+.$$
 It follows from (\ref{cor2}) and Lemma \ref{PIPI} that
 \BQN\label{uniform}
 \lim_{u\rw\IF}\sup_{(k,l)\in K_u^+}\left|\frac{\pk{\sup_{(s,t)\in I_{0,0}(u)}X_{u,k,l}(s,t)>u_{k,l}^{-\epsilon}}}{\Psi(u_{k,l}^{-\epsilon})}-\prod_{i=1}^2\mathcal{H}_{\alpha_i}[0,S]\right|=0.
 \EQN
 Further,
 \BQN\label{medium}
 \sum_{(k,l)\in K_u^+} \pk{\sup_{(s,t)\in I_{k,l}(u)}\overline{X}(s,t)>u_{k,l}^{-\epsilon}}&\sim& \prod_{i=1}^2\mathcal{H}_{\alpha_i}[0,S]\sum_{(k,l)\in K_u^+} \Psi(u_{k,l}^{-\epsilon})\nonumber\\
 &\sim& \prod_{i=1}^2\mathcal{H}_{\alpha_i}[0,S]\Psi(u)\sum_{(k,l)\in K_u^+} e^{-(1-\epsilon)u^2\inf_{s\in I_{k,l}(u)}v^2(|s+bt|)}\nonumber\\
 &\sim& \prod_{i=1}^2\mathcal{H}_{\alpha_i}[0,S]\Psi(u)\sum_{|l|\leq N_2(u)+1}\sum_{k\in E_l^+} e^{-(1-\epsilon)u^2\inf_{(s,t)\in I_{k,l}(u)}v^2(|s+bt|)}.
 \EQN
 Note that
 $$\{k:|k+a_l(u)|\leq N_1(u)-a_1(u)-2\}\subset E_l^-\subset E_l^+\subset \{k:|k+a_l(u)|\leq N_1(u)+a_{1}(u)+2\},$$
 with $a_l(u)=\left[\frac{b\la{\rho}_2(u^{-1})l}{\la{\rho}_1(u^{-1})}\right], l\in \mathbb{R}.$

By $\eta\in (0,\IF)$, we have
$$\lim_{u\rw\IF}\frac{\la{\rho}_2(u^{-1})}{\la{\rho}_1(u^{-1})}=\eta^{-1/\alpha_1},$$
which implies that for $2|b|\eta^{-1/\alpha}+2\leq L\leq |k+a_l(u)|\leq N_1(u)+a_{1}(u)+2$ and $u$ large enough,
$$(L-2|b|\eta^{-1/\alpha}-2)\la{\rho}_1(u^{-1}) S\leq \inf_{(s,t)\in I_{k,l}(u)}|s+bt|\leq \sup_{(s,t)\in I_{k,l}(u)}|s+bt|\leq (L+2|b|\eta^{-1/\alpha}+2)\la{\rho}_1(u^{-1}) S .$$
By Lemma 6.1 in \cite{KEP20162}, we have that for any $0<\epsilon<1$, as $u$ sufficiently large,
\BQNY
\frac{v^2(|s+bt|)}{v^2(|s'+bt'|)}&\geq& (1-\epsilon/2)\min\left(\left|\frac{s+bt}{s'+bt'}\right|^{\beta-\epsilon}, \left|\frac{s+bt}{s'+bt'}\right|^{\beta+\epsilon}\right)\\
&\geq&(1-\epsilon/2)\left(\frac{L-2|b|\eta^{-1/\alpha}-2}{L+2|b|\eta^{-1/\alpha}+2}\right)^{\beta+\epsilon}, \quad (s,t), (s',t')\in I_{k,l}(u)
\EQNY
with $\quad L\leq |k+a_l(u)|\leq N_1(u)+a_{1}(u)+2, |l|\leq N_2(u)+1$.
Thus for any $0<\epsilon<1$ there exists $k_\epsilon>1$ such that for $u$ large enough,
\BQN\label{mod1}\inf_{s\in I_{k,l}(u)}v^2(|s+bt|)\geq (1-\epsilon)\sup_{s\in I_{k,l}(u)}v^2(|s+bt|), \quad k_\epsilon\leq |k+a_l(u)|\leq N_1(u)+a_{1}(u)+2, |l|\leq N_2(u)+1.
\EQN
Consequently, by  Lemma 6.3 in \cite{KEP20162} and as a continuation of (\ref{medium}) we have that
 \BQN\label{upper}
 &&\sum_{(k,l)\in K_u^+} \pk{\sup_{(s,t)\in I_{k,l}(u)}\overline{X}(s,t)>u_{k,l}^{-\epsilon}}\nonumber\\
 && \ \ \leq \prod_{i=1}^2\mathcal{H}_{\alpha_i}[0,S]\Psi(u)\sum_{|l|\leq N_2(u)+1}\left(2k_\epsilon+1\right.\nonumber\\
 && \ \ \ \ \left.+\frac{1}{S^2\la{\rho}_1(u^{-1})\la{\rho}_2(u^{-1})}\sum_{k_\epsilon\leq |k+a_l(u)|\leq N_1(u)+a_1(u)+2}\int_{(s,t)\in I_{k,l}(u)}e^{-(1-\epsilon)^2u^2v^2(|s+bt|)}dsdt\right)\nonumber\\
 && \ \ \leq \prod_{i=1}^2\mathcal{H}_{\alpha_i}[0,S]\Psi(u)\sum_{|l|\leq N_2(u)+1}\left(2k_\epsilon+1+\frac{2}{S\la{\rho}_1(u^{-1})}\int_{0}^{\la{v}(\ln u/u)}e^{-(1-\epsilon)^2u^2v^2(|s|)}ds\right)\nonumber\\
 && \ \ \sim \prod_{i=1}^2\mathcal{H}_{\alpha_i}[0,S]\Psi(u)2(N_2(u)+1)\left(2k_\epsilon+1+\frac{2}{S\la{\rho}_1(u^{-1})}\int_{0}^{\la{v}(\ln u/u)}e^{-(1-\epsilon)^2u^2v^2(|s|)}ds\right)\nonumber\\
 && \ \ \sim \prod_{i=1}^2\mathcal{H}_{\alpha_i}[0,S]\Psi(u)2(N_2(u)+1)\left(2k_\epsilon+1+\frac{2(1-\epsilon)^{-2/\beta}\Gamma(1/\beta+1)
 \la{v}(u^{-1})}{S\la{\rho}_1(u^{-1})}\right)\nonumber\\
 && \ \ \sim \prod_{i=1}^2\mathcal{H}_{\alpha_i}\frac{4T_2\Gamma(1/\beta+1)\la{v}(u^{-1})}{\la{\rho}_1(u^{-1})\la{\rho}_2(u^{-1})}\Psi(u), \quad u\rw\IF, S\rw\IF, \epsilon\rw 0.
 \EQN
 Similarly,
 \BQN\label{lower}
 \sum_{(k,l)\in K_u^-} \pk{\sup_{(s,t)\in I_{k,l}(u)}\overline{X}(s,t)>u_{k,l}^{+\epsilon}}\sim \prod_{i=1}^2\mathcal{H}_{\alpha_i}\frac{4T_2\Gamma(1/\beta+1)\la{v}(u^{-1})}{\la{\rho}_1(u^{-1})\la{\rho}_2(u^{-1})}\Psi(u), \quad u\rw\IF, S\rw\IF, \epsilon\rw 0.
 \EQN
 Next we will show that $\Sigma_i(u), i=1,2,3,$ are all negligible in comparison with $\sum_{(k,l)\in K_u^+} \pk{\sup_{(s,t)\in I_{k,l}(u)}\overline{X}(s,t)>u_{k,l}^{-\epsilon}}$. For any $(k,l,k_1,l_1)\in \Lambda_1$, without loss of generality, we assume that $k+1=k_1$. Let
 $$ I_{k,l}^1(u)=[k\la{\rho}_1(u^{-1})S, \la{\rho}_1(u^{-1})((k+1)S-\sqrt{S})]\times J_l(u), \quad I_{k,l}^2=[ \la{\rho}_1(u^{-1})((k+1)S-\sqrt{S}), (k+1)\la{\rho}_1(u^{-1})S]\times J_l(u),$$
 then
 $I_{k,l}(u)=I_{k,l}^{1}(u)\cup I_{k,l}^2(u)$ and
\BQNY
&&\pk{\sup_{(s,t)\in I_{k,l}(u)}\overline{X}(s,t)>u_{k,l}^{-\epsilon}, \sup_{(s,t)\in I_{k_1,l_1}(u)}\overline{X}(s,t)>u_{k_1,l_1}^{-\epsilon}}\\
&& \  \ \leq
\pk{\sup_{(s,t)\in I_{k,l}^1(u)}\overline{X}(s,t)>u_{k,l}^{-\epsilon}, \sup_{(s,t)\in I_{k_1,l_1}(u)}\overline{X}(s,t)>u_{k_1,l_1}^{-\epsilon}}+\pk{\sup_{(s,t)\in I_{k,l}^2(u)}\overline{X}(s,t)>u_{k,l}^{-\epsilon}}
\EQNY
Analogously  as in (\ref{uniform}), we have
\BQNY
 \lim_{u\rw\IF}\sup_{(k,l)\in K_u^-}\left|\frac{\pk{\sup_{(s,t)\in I_{k,l}^2(u)}\overline{X}(s,t)>u_{k,l}^{-\epsilon}}}{\Psi(u_{k,l}^{-\epsilon})}-\mathcal{H}_{\alpha_1}[0,\sqrt{S}]\mathcal{H}_{\alpha_2}[0,S]\right|=0.
 \EQNY

Moreover£¬, in light of (\ref{cor2}) and Lemma 5.4 in \cite{KEP20162} (or Corollary 3.2 in \cite{KEP20161}) we have for $u$ large enough,
$$\pk{\sup_{(s,t)\in I_{k,l}^1(u)}\overline{X}(s,t)>u_{k,l}^{-\epsilon}, \sup_{(s,t)\in I_{k_1,l_1}(u)}\overline{X}(s,t)>u_{k_1,l_1}^{-\epsilon}}\leq \mathcal{C}S^4e^{-\mathcal{C}_1S^{\alpha^*/4}}\Psi(u_{k,l,k_1,l_1}^{-\epsilon}),$$
and for $(k,l,k_1,l_1)\in \Lambda_2$,
$$\pk{\sup_{(s,t)\in I_{k,l}(u)}\overline{X}(s,t)>u_{k,l}^{-\epsilon}, \sup_{(s,t)\in I_{k_1,l_1}(u)}\overline{X}(s,t)>u_{k_1,l_1}^{-\epsilon}}\leq \mathcal{C}S^4e^{-\mathcal{C}_1|(k-k_1)^2+(l-l_1)^2|^{\alpha^*/4}S^{\alpha^*/2}}\Psi(u_{k,l,k_1,l_1}^{-\epsilon}),$$
with $u_{k,l,k_1,l_1}^{-\epsilon}=\min (u_{k,l}^{-\epsilon}, u_{k_1,l_1}^{-\epsilon})$, $\alpha^*=\min(\alpha_1,\alpha_2)$ and $\mathcal{C}$, $\mathcal{C}_1$ being some fixed positive constants independent of $S$ and $u$.
Since each $I_{k,l}(u)$ has at most $8$ neighbors, then
\BQN\label{double1}
\Sigma_1(u)&\leq 2&\sum_{(k,l,k_1,l_1)\in \Lambda_1(u)}\left(\mathcal{H}_{\alpha_1}[0,\sqrt{S}]\mathcal{H}_{\alpha_2}[0,S]+\mathcal{H}_{\alpha_1}[0,S]\mathcal{H}_{\alpha_2}[0,\sqrt{S}]
+\mathcal{C}S^4e^{-\mathcal{C}_1S^{\alpha^*/4}}\right)\Psi(u_{k,l,k_1,l_1}^{-\epsilon})\nonumber\\
&\leq& 16\sum_{(k,l)\in K_u^-}\left(\mathcal{H}_{\alpha_1}[0,\sqrt{S}]\mathcal{H}_{\alpha_2}[0,S]+\mathcal{H}_{\alpha_1}[0,S]\mathcal{H}_{\alpha_2}[0,\sqrt{S}]
+\mathcal{C}S^4e^{-\mathcal{C}_1S^{\alpha^*/4}}\right)\Psi(u_{k,l}^{-\epsilon})\nonumber\\
&=&o\left( \frac{\la{v}(u^{-1})}{\la{\rho}_1(u^{-1})\la{\rho}_2(u^{-1})}\Psi(u)\right), \ \ u\rw\IF, S\rw\IF,
\EQN
and
\BQN\label{double2}
\Sigma_2(u)&\leq&\sum_{(k,l,k_1,l_1)\in \Lambda_2(u)}\mathcal{C}S^4e^{-\mathcal{C}_1|(k-k_1)^2+(l-1_1)^2|^{\alpha^*/4}S^{\alpha^*/2}}\Psi(u_{k,l,k_1,l_1}^{-\epsilon})\nonumber\\
&\leq&\sum_{(k,l)\in K_u^-}\Psi(u_{k,l}^{-\epsilon})\mathcal{C}S^4\sum_{|k_1|+|l_1|\geq 1}e^{-\mathcal{C}_1|k_1^2+l_1^2|^{\alpha^*/4}S^{\alpha^*/2}}\nonumber\\
&\leq&\sum_{(k,l)\in K_u^-}\Psi(u_{k,l}^{-\epsilon})\mathcal{C}S^4e^{-\mathbb{Q}S^{\alpha^*/2}}=o\left( \frac{\la{v}(u^{-1})}{\la{\rho}_1(u^{-1})\la{\rho}_2(u^{-1})}\Psi(u)\right), \ \ u\rw\IF, S\rw\IF.
\EQN
For $(k,l,k_1,l_1)\in \Lambda_3(u)$, $|t-t_1|\geq \epsilon/2$ holds with $(s,t)\in I_{k,l}(u)$ and $(s_1,t_1)\in I_{k_1,l_1}(u)$. Then by (\ref{local}), for $u$ large enough,
$$\Var(\overline{X}(s,t)+\overline{X}(s_1,t_1))=2(1+r(s,t,s_1,t_1))\leq 2+2\sup_{|t-t_1|\geq \epsilon/2}r(s,t,s_1,t_1)\leq 4-\delta$$
holds with $0<\delta<1$ for $(k,l,k_1,l_1)\in \Sigma_3(u),(s,t)\in I_{k,l}(u), (s_1,t_1)\in I_{k_1,l_1}(u)$. Further, Borell-TIS inequality leads to
\BQN\label{double3}
\Sigma_3(u)&\leq& \sum_{(k,l,k_1,l_1)\in \Lambda_3(u)}\pk{\sup_{(s,t,s_1,t_1)\in I_{k,l}(u)\times I_{k_1,l_1}(u)}\overline{X}(s,t)+\overline{X}(s_1,t_1)>2u}\nonumber\\
&\leq&\sum_{(k,l,k_1,l_1)\in \Lambda_3(u)}e^{-\frac{(2u-a_1)^2}{2(4-\delta)}}
\leq \frac{\mathbb{Q}}{\la{\rho}_1(u^{-1})\la{\rho}_2(u^{-1})}e^{-\frac{(2u-a_1)^2}{2(4-\delta)}}=o\left( \frac{\la{v}(u^{-1})}{\la{\rho}_1(u^{-1})\la{\rho}_2(u^{-1})}\Psi(u)\right), \quad u\rw\IF,
\EQN
with $a_1=2\mathbb{E}\left(\sup_{(s,t)\in E\cap \{(s,t):|s+bt|\leq \epsilon\}}\overline{X}(s,t)\right)$. Inserting (\ref{upper})-(\ref{double3}) into
(\ref{main1}) yields that
$$\pi_1(u)\sim 4T_2\Gamma(1/\beta+1)\prod_{i=1}^2\mathcal{H}_{\alpha_i}\frac{\la{v}(u^{-1})}{\la{\rho}_1(u^{-1})\la{\rho}_2(u^{-1})}\Psi(u), \quad u\rw\IF.$$
This establishes the claim.\\
 \underline{$\gamma_1\in (0,\IF)$}. We first introduce some new notation for further analysis. Let
\BQNY
\widehat{J}_{l}(u)&=&\{|s+bt|\leq \la{\rho}_1(u^{-1})S, l\la{\rho}_2(u^{-1})S\leq t\leq (l+1)\la{\rho}_2(u^{-1})S\},\\
  J_{k,l}(u)&=&\{k \la{\rho}_1(u^{-1})S\leq s+bt\leq (k+1)\la{\rho}_1(u^{-1})S, l\la{\rho}_2(u^{-1})S\leq t\leq (l+1)\la{\rho}_2(u^{-1})S\},\\
  L_u&=&\{(k,l):  k\neq 0, -1, |k|\leq N_1(u)+1, |l|\leq N_2(u)+1\}, \quad u_{k}^{-\epsilon}=u\left(1+(1-\epsilon)\inf_{s\in I_{k}(u)}v^2(|s|)\right).
  \EQNY
  Using Bonferroni inequality, we have
  \BQN
  \pi_2(u)-\sum_{i=4}^5\Sigma_i(u)\leq \pi_1(u)\leq \pi_2(u)+\sum_{(k,l)\in L_u}\pk{\sup_{(s,t)\in J_{k,l}(u)}\overline{X}(s,t)>u_{k}^{-\epsilon}},
  \EQN
 where
 $$\pi_2(u)=\sum_{|l|\leq N_2(u)+1}\pk{\sup_{(s,t)\in \widehat{J}_{l}(u)}\frac{\overline{X}(s,t)}{1+(1+o(1))v^2(|s+bt|)}>u},$$
 $$\Sigma_4(u)=\sum_{|l|\leq N_2(u)+1}\pk{\sup_{(s,t)\in \widehat{J}_{l}(u)}\overline{X}(s,t)>u, \sup_{(s,t)\in \widehat{J}_{l+1}(u)}\overline{X}(s,t)>u},$$
  $$\Sigma_5(u)=\sum_{|l|, |l_1|\leq N_2(u)+1, l+2\leq l_1}\pk{\sup_{(s,t)\in \widehat{J}_{l}(u)}\overline{X}(s,t)>u, \sup_{(s,t)\in \widehat{J}_{l_1}(u)}\overline{X}(s,t)>u}.$$
  In order to apply Lemma \ref{PIPI}, we set
  $$X_{u,l}(s,t)=\overline{X}(s-bl\la{\rho}_2(u^{-1})S,l\la{\rho}_2(u^{-1})S+t), (s,t)\in  \widehat{J}_{0}(u), |l|\leq N_2(u)+1, $$
  By (\ref{cor2}) and the fact that
  $$ u^2 v^2(|\la{\rho}_1(u^{-1})s+b\la{\rho}_2(u^{-1})t|)\rw \gamma_1|s+b\eta^{-1/\alpha_1}t|^\beta, \quad u\rw\IF,$$
  holds uniformly with respect to $(s,t)\in \{(s,t): |s+b\eta^{-1/\alpha_1}t|\leq S, 0\leq t\leq S\}$, and using Lemma \ref{PIPI}, we have
  \BQNY
 \lim_{u\rw\IF}\sup_{|l|\leq N_2(u)+1}\left|\frac{\pk{\sup_{(s,t)\in \widehat{J}_{0}(u)}\frac{X_{u,l}(s,t)}{1+(1+o(1))v^2(|s+bt|)}>u}}{\Psi(u)}-\mathcal{H}_{\alpha_1}^{\gamma_1, b\eta^{-1/\alpha_1}}(S)\right|=0,
 \EQNY
 where
$\mathcal{H}_{\alpha_1}^{\gamma_1, b\eta^{-1/\alpha_1}}(S)$ is defined in (\ref{piterbarg}).
 Thus
 \BQN\label{eq33}
 \pi_2(u)\sim \sum_{|l|\leq N_2(u)+1}\mathcal{H}_{\alpha_1}^{\gamma_1, b\eta^{-1/\alpha_1}}(S)\Psi(u)\sim \mathcal{H}_{\alpha_1}^{\gamma_1, b\eta^{-1/\alpha_1}}\frac{2T_2}{\la{\rho}_2(u^{-1})}\Psi(u), \quad u\rw\IF.
 \EQN

 Observe that for $u$ sufficiently large,
 \BQN\label{mod3}
 J_{k,l}(u)&\subset& \{(k-\mu) \la{\rho}_1(u^{-1})S-bl \la{\rho}_2(u^{-1})S \leq s\leq (k+1+\mu) \la{\rho}_1(u^{-1})S-bl \la{\rho}_2(u^{-1})S\}\\
 && \times J_l(u):=\widehat{I}_{k,l}(u).\nonumber
 \EQN
 with $\mu=2|b|\eta^{-1/\alpha_1}$.
 Let
  $$X_{u,k,l}(s,t)=\overline{X}(s+k\la{\rho}_1(u^{-1})S-bl\la{\rho}_2(u^{-1})S,l\la{\rho}_2(u^{-1})S+t), (s,t)\in  \widehat{I}_{0,0}(u), (k,l)\in L_u.$$
  Thus in light of (\ref{cor2}) and Lemma \ref{PIPI} we have
  \BQNY
 \lim_{u\rw\IF}\sup_{(k,l)\in L_u}\left|\frac{\pk{\sup_{(s,t)\in \widehat{I}_{0,0}(u)}X_{u,k,l}(s,t)>u_{k}^{-\epsilon}}}{\Psi(u_k^{-\epsilon})}-\mathcal{H}_{\alpha_1}[-\mu S, (\mu+1) S]\mathcal{H}_{\alpha_2}[0,S]\right|=0.
 \EQNY
Moreover, Potter's bound (see, e.g., \cite{BI1989}) shows that for $u$ large enough and $S>1$
 $$u^2v^2(|s|)\geq\frac{\gamma_1}{2}\frac{\rho_1^2(\la{\rho}_1(u^{-1})|\frac{s}{\la{\rho}_1(u^{-1})}|)}{\rho_1^2(\la{\rho}_1(u^{-1}))}\geq \mathbb{Q}|kS|^{\beta/2}, \quad s\in I_k(u), |k|\leq N_1(u)+1, k\neq -1,0 .$$
 Consequently,
 \BQN\label{upper1}
 &&\sum_{(k,l)\in L_u}\pk{\sup_{(s,t)\in J_{k,l}(u)}\overline{X}(s,t)>u_{k}^{-\epsilon}}\nonumber\\
 && \ \ \leq \sum_{(k,l)\in L_u}\pk{\sup_{(s,t)\in \widehat{I}_{k,l}(u)}\overline{X}(s,t)>u_{k}^{-\epsilon}}\nonumber\\
 && \ \ = \sum_{(k,l)\in L_u}\pk{\sup_{(s,t)\in \widehat{I}_{0,0}(u)}X_{u,k,l}(s,t)>u_{k}^{-\epsilon}}\nonumber\\
 && \ \  \sim \sum_{(k,l)\in L_u}\mathcal{H}_{\alpha_1}[-\mu S, (\mu+1) S]\mathcal{H}_{\alpha_2}[0,S]\Psi(u_k^{-\epsilon})\nonumber\\
 && \ \  \sim \mathcal{H}_{\alpha_1}[-\mu S, (\mu+1) S]\mathcal{H}_{\alpha_2}[0,S]\frac{2T_2}{S\la{\rho}_2(u^{-1})}\Psi(u)\sum_{|k|\leq N_1(u)+1, k\neq -1,0} e^{-(1-\epsilon)u^2\inf_{s\in I_k(u)}v^2(|s|)}\nonumber\\
 && \ \ \leq \mathcal{H}_{\alpha_1}[-\mu S, (\mu+1) S]\mathcal{H}_{\alpha_2}[0,S]\frac{2T_2}{S\la{\rho}_2(u^{-1})}\Psi(u) \sum_{|k|\leq N_1(u)+1, k\neq -1,0} e^{-\mathbb{Q}|kS|^{\beta/2}}\nonumber\\
 && \ \ \leq \mathcal{H}_{\alpha_1}[-\mu S, (\mu+1) S]\mathcal{H}_{\alpha_2}[0,S]\frac{2T_2}{S\la{\rho}_2(u^{-1})}\Psi(u) e^{-\mathbb{Q}_1S^{-\beta/2}}=o\left(\pi_2(u)\right), \quad u\rw\IF, S\rw \IF.
 \EQN
 Analogously as in (\ref{double1})-(\ref{double3}), we get that
 $$\Sigma_{i}(u)=o(\pi_2(u)), i=4,5, \quad u\rw\IF, S\rw\IF.$$
 Therefore, we conclude that
\BQN\label{conc}
\pi_1(u)\sim   \mathcal{H}_{\alpha_1}^{\gamma_1, b\eta^{-1/\alpha_1}}\frac{2T_2}{\la{\rho}_2(u^{-1})}\Psi(u), \quad u\rw\IF.
\EQN
 This establishes the claim.\\
 \underline{$\gamma_1=\IF$}. For any $y>0$, we have for $u$ sufficiently large,
 \BQN
 \pk{\sup_{t\in [-T_2,T_2]}X(-b t,t)>u}\leq\pi_1(u)\leq \pi_3^y(u)+\sum_{(k,l)\in L_u}\pk{\sup_{(s,t)\in J_{k,l}(u)}\overline{X}(s,t)>u_{k}^{-\epsilon}},
  \EQN
  where
  $$\pi_3^y(u)=\sum_{|l|\leq N_2(u)+1}\pk{\sup_{(s,t)\in \widehat{J}_{l}(u)}\frac{\overline{X}(s,t)}{1+y\rho_1^2(|s+bt|)}>u}.$$
 By (\ref{eq33}), we have
 \BQNY
 \pi_3^y(u)\sim  \frac{\mathcal{H}_{\alpha_1}^{y, b\eta^{-1/\alpha_1}}(S)}{S}\frac{2T_2}{\la{\rho}_2(u^{-1})}\Psi(u), \ \ u\rw\IF.
 \EQNY
Duo to
 $$\lim_{y\rw\IF}\mathcal{H}_{\alpha_1}^{y, b\eta^{-1/\alpha_1}}(S)=\mathbb{E}\left(e^{\sup_{t\in[0,S]}W_{\alpha_1, \alpha_1}(-b\eta^{-1/\alpha_1}t,t)}\right)=\mathcal{H}_{\alpha_1}
 [0,(|b|^{\alpha_1}\eta^{-1}+1)^{1/\alpha_1}S],$$
 we have
 \BQN\label{eq41}
 \pi_3^y(u)\sim  (|b|^{\alpha_1}\eta^{-1}+1)^{1/\alpha_1}\mathcal{H}_{\alpha_1}\frac{2T_2}{\la{\rho}_2(u^{-1})}\Psi(u), \ \ u\rw\IF, y\rw\IF, S\rw\IF.
 \EQN
 Observe that
 $X(-b t,t), t\in [-T_2,T_2]$ is a Gaussian process with variance $1$ and correlation function satisfying
 $$1-Cor(X(-b t,t), X(-b s,s))\sim \rho_1^2(|b(t-s)|)+\rho_2^2(|t-s|)\sim (|b|^{\alpha_1}\eta^{-1}+1)\rho_2^2(|t-s|), \quad |t-s|\rw 0,$$
  and $$Cor(X(-b t,t), X(-b s,s))<1, \quad t\neq s, s,t \in [-T_2,T_2].$$
  In light of Lemma 7.1 in \cite{Pit96} and substituting the polynomial function by $(|b|^{\alpha_1}\eta^{-1}+1)^{1/\alpha_1}/\la{\rho}_2(u^{-1})$, we have
  $$\pk{\sup_{t\in [-T_2,T_2]}X(-b t,t)}\sim (|b|^{\alpha_1}\eta^{-1}+1)^{1/\alpha_1}\mathcal{H}_{\alpha_1}\frac{2T_2}{\la{\rho}_2(u^{-1})}\Psi(u),$$
  which combined with (\ref{upper1}) and (\ref{eq41}) completes the proof.  \QED

  \prooftheo{Th1}. The proof of this theorem follows line by line the same as  the proof of Theorem \ref{Th2} by letting $b=0$ except some tiny modification. Here we just point out several arguments that require special attention.\\
  \underline{$\gamma_1=0$}.
 The derivation of (\ref{mod1}):
 By Lemma 6.1 in \cite{KEP20162}, we have that for any $0<\epsilon<1$, as $u$ sufficiently large,
\BQNY
\frac{v^2(|s|)}{v^2(|s'|)}&\geq& (1-\epsilon/2)\min\left(\left|\frac{s}{s'}\right|^{\beta-\epsilon}, \left|\frac{s}{s'}\right|^{\beta+\epsilon}\right)\\
&\geq&(1-\epsilon/2)\left(\frac{|k|-1}{|k|+1}\right)^{\beta+\epsilon}, \quad s,s'\in I_{k}(u), 2\leq |k| \leq N_1(u)+2.
\EQNY
 Thus for any $0<\epsilon<1$, there exists $k_\epsilon>2$ such that for $u$ large enough,
 $$\inf_{s\in I_k(u)}v^2(|s|)\geq (1-\epsilon)\sup_{s\in I_k(u)}v^2(|s|), \quad k_\epsilon\leq |k|\leq N_1(u)+2.$$
 \underline{$\gamma_1\in (0,\IF)$}. Note that from (\ref{eq33}) to (\ref{conc}),  $b\eta^{-1/\alpha_1}\big{|}_{b=0}$ may have no meaning since $\lim_{t\downarrow 0}\frac{\rho_2^2(t)}{\rho_1^2(t)}\in [0,\IF]$. Replacing  $b\eta^{-1/\alpha_1}$ by $0$ from (\ref{eq33}) to (\ref{conc}), we have
 $$\mathcal{H}_{\alpha_1}^{\gamma_1, b\eta^{-1/\alpha_1}}(S)=\mathcal{P}_{\alpha_1}^{\gamma_1}[-S,S]\mathcal{H}_{\alpha_2}[0,S], \quad \mathcal{H}_{\alpha_1}^{\gamma_1, b\eta^{-1/\alpha_1}}=\mathcal{P}_{\alpha_1}^{\gamma_1}\mathcal{H}_{\alpha_2}.$$
 Moreover, $\mu$, which appears from (\ref{mod3}) to (\ref{upper1}),  should also be substituted by $0$.\\
 \underline{$\gamma_1=\IF$}. Letting $b=0$ and $b\eta^{-1/\alpha_1}=0$ in the proof of $\gamma_1=\IF$ in Theorem \ref{Th2} establishes the claim. This completes the proof. \QED

\prooftheo{th31} Due to {\bf F}, we know that $f'(t)\neq 0$ for $t\in (T_1,T_2)$.   Without loss of generality, we assume that $f'(t)>0, t\in (T_1, T_2)$.
 Let $Y(s,t)=X(s,\la{f}(s-t)), (s,t)\in E_1:=\{(s,t): S_1\leq s\leq S_2, f(T_1)\leq s-t\leq f(T_2)\}$. Then
\BQNY
\pk{\sup_{(s,t)\in [S_1,S_2]\times [T_1,T_2] }X(s,t)>u}=\pk{\sup_{(s,t)\in E_1}Y(s,t)>u}.
\EQNY
Clearly, $\sigma_Y(s,t), (s,t)\in E_1$ attains its maximum over the line $\{(s,0): f(T_1)\leq s\leq f(T_2)\}$. Moreover, by (\ref{cor4}), (\ref{var4}) and Lemma 6.4 in \cite{KEP20162}, we have
\BQNY
1-r_Y(s,t,s',t')\sim \rho_1^2(|s-s_1|)+\rho_2^2\left(\frac{|s-s'-(t-t')|}{f'(\la{f}(s))}\right)\sim  \rho_1^2(|s-s'|)+\rho_2^2\left(\frac{|t-t'|}{f'(\la{f}(s))}\right), \quad t,t'\rw 0, |s-s'|\rw 0,
\EQNY
and
\BQNY
1-\sigma_Y(s,t)\sim v^2(g(\la{f}(s))|t|), \quad |t|\rw 0.
\EQNY
Let $E_2=[f(T_1)+\delta, f(T_2)-\delta]\times [-\delta, \delta]$, $E_3=[f(T_1)-\delta, f(T_1)+\delta]\times [-\delta, \delta]$ and $E_4=[ f(T_2)-\delta, f(T_2)+\delta]\times [-\delta, \delta]$. Then for $u$ large enough, we have
\BQN\label{main4}
\pi_{4,\delta}(u)\leq \pk{\sup_{(s,t)\in E_1}Y(s,t)>u}\leq \sum_{i=4}^6\pi_{i,\delta}(u)+\pk{\sup_{(s,t)\in E_1\setminus\{(s,t): |t|\leq \delta\}}Y(s,t)>u}
\EQN
where
$$\pi_{i,\delta}(u)=\pk{\sup_{(s,t)\in E_1\cap E_{i-2}}Y(s,t)>u}, \quad i=4,5,6.$$
We first focus on $\pi_{4,\delta}(u)$. Denote by
$F_k=\left[s_k, s_{k+1}\right],   0\leq k\leq n $
with $s_k=f(T_1)+\delta+\frac{k(f(T_2)-f(T_1)-2\delta)}{n}$.
Then we have
\BQN\label{eq6}
\sum_{k=0}^{n}\pk{\sup_{(s,t)\in F_k\times [-\delta, \delta]}Y(s,t)>u}-\sum_{i=6}^7\Sigma_i(u)\leq \pi_{4,\delta}(u)\leq \sum_{k=0}^{n}\pk{\sup_{(s,t)\in F_k\times [-\delta, \delta]}Y(s,t)>u},
\EQN
where
\BQNY
\Sigma_6(u)&=&\sum_{0\leq k\leq n-1}\pk{\sup_{(s,t)\in F_k\times [-\delta, \delta]}Y(s,t)>u, \sup_{(s,t)\in F_{k+1}\times [-\delta, \delta]}Y(s,t)>u }, \\
 \Sigma_7(u)&=&\sum_{0\leq k+1<l\leq n}\pk{\sup_{(s,t)\in F_k\times [-\delta, \delta]}Y(s,t)>u, \sup_{(s,t)\in F_l\times [-\delta, \delta]}Y(s,t)>u }.
\EQNY
Let $Z(s,t)$ be a homogeneous Gaussian field with variance $1$ and correlation function satisfying
\BQN\label{corr}
1-r_Z(s,t,s',t')\sim \rho_1^2(|s-s'|)+\rho_2^2\left(|t-t'|\right), \quad |s-s'|, |t-t'|\rw 0.
\EQN
 For any $0<\epsilon<1/2$, if $n$ sufficiently large and $\delta$ sufficiently small,
\BQNY
&&1-r_Z\left((1-\epsilon)s,(1-\epsilon)(f'(\la{f}(s_k)))^{-1}t,(1-\epsilon)s',(1-\epsilon)(f'(\la{f}(s_k)))^{-1}t'\right)\\
&& \quad\leq 1-r_Y(s,t,s',t')\\
&& \quad \leq 1-r_Z\left((1+\epsilon)s,(1+\epsilon)(f'(\la{f}(s_k)))^{-1}t,(1+\epsilon)s',(1+\epsilon)(f'(\la{f}(s_k)))^{-1}t'\right),
\EQNY
and
\BQNY
(1-\epsilon)|g(\la{f}(s_k))|^\beta v^2(|t|)
 \leq 1-\sigma_Y(s,t)\leq (1+\epsilon)|g(\la{f}(s_k))|^\beta v^2(|t|)
\EQNY
hold for $(s,t), (s',t')\in F_k\times [-\delta, \delta]$ and $1\leq k\leq n$. Thus by Slepian inequality (\cite{AdlerTaylor}), we have
\BQNY
\pk{\sup_{(s,t)\in F_k\times [-\delta, \delta]}\frac{Z((1-\epsilon)s,(1-\epsilon)(f'(\la{f}(s_k)))^{-1}t)}{1+(1+\epsilon)|g(\la{f}(s_k))|^\beta v^2(|t|)}>u}&\leq& \pk{\sup_{(s,t)\in F_k\times [-\delta, \delta]}Y(s,t)>u}\\
&\leq& \pk{\sup_{(s,t)\in F_k\times [-\delta, \delta]}\frac{Z((1+\epsilon)s,(1+\epsilon)(f'(\la{f}(s_k)))^{-1}t)}{1+(1-\epsilon)|g(\la{f}(s_k))|^\beta v^2(|t|)}>u}.
\EQNY
Direct calculation yields that
\BQNY
&&1-r_Z\left((1+\epsilon)s,(1+\epsilon)(f'(\la{f}(s_k)))^{-1}t,(1+\epsilon)s',(1+\epsilon)(f'(\la{f}(s_k)))^{-1}t'\right)\\
&&\quad \sim(1+\epsilon)^{\alpha_1}\rho_1^2(|s-s'|)+(1+\epsilon)^{\alpha_2}(f'(\la{f}(s_k)))^{-\alpha_2}\rho_2^2\left(|t-t'|\right), \quad |s-s'|, |t-t'|\rw 0.
\EQNY
{\it \underline{Case $\gamma_2=0$}}.
For simplicity, we denote by
\BQN\label{theta}
\Theta(u):=2\Gamma(1/\beta+1)\left(\prod_{i=1}^2\mathcal{H}_{\alpha_i}\right)\frac{\overleftarrow{\vv}(1/u)}
{\overleftarrow{\LL}_1(1/u)\overleftarrow{\LL}_2(1/u)}\Psi(u).
\EQN
By  Theorem \ref{Th1}, we have that
\BQNY
\pk{\sup_{(s,t)\in F_k\times [-\delta, \delta]}\frac{Z((1+\epsilon)s,(1+\epsilon)(f'(\la{f}(s_k)))^{-1}t)}{1+(1-\epsilon)|g(\la{f}(s_k))|^\beta v^2(|t|)}>u}
 \sim a_1(\epsilon)\frac{s_{k+1}-s_k}{|g(\la{f}(s_k))f'(\la{f}(s_k))|}
\Theta(u),
\EQNY
 with $a_1(\epsilon)=(1+\epsilon)^2(1-\epsilon)^{-1/\beta}$.
Further, for $u\rw\IF, n\rw\IF$,
\BQNY
&&\sum_{k=0}^{n}\pk{\sup_{(s,t)\in F_k\times [-\delta, \delta]}Y(s,t)>u}\\
&& \quad \leq a_1(\epsilon)\left(\sum_{k=0}^{n} \frac{s_{k+1}-s_k}{|g(\la{f}(s_k))f'(\la{f}(s_k))|} \right)
\Theta(u)(1+o(1))\\
&& \quad \leq a_1(\epsilon)\int_{f(T_1)+\delta}^{f(T_2)-\delta}\left|g(\la{f}(s))f'(\la{f}(s))\right|^{-1}ds
\Theta(u)(1+o(1)).
\EQNY
Similarly, as $u\rw\IF, n\rw\IF$,
\BQNY
&&\sum_{k=0}^{n}\pk{\sup_{(s,t)\in F_k\times [-\delta, \delta]}Y(s,t)>u}\\
&& \quad \geq a_1(-\epsilon)\int_{f(T_1)+\delta}^{f(T_2)-\delta}\left|g(\la{f}(s))f'(\la{f}(s))\right|^{-1}ds
\Theta(u)(1+o(1)).
\EQNY
Next we focus on $\Sigma_6(u)$ and $\Sigma_7(u)$. By the fact that
\BQN\label{eq3}
&&\pk{\sup_{(s,t)\in F_k\times [-\delta, \delta]}Y(s,t)>u, \sup_{(s,t)\in F_{k+1}\times [-\delta, \delta]}Y(s,t)>u }\nonumber\\
&&\quad \leq \pk{\sup_{(s,t)\in F_k\times [-\delta, \delta]}Y(s,t)>u}+\pk{ \sup_{(s,t)\in F_{k+1}\times [-\delta, \delta]}Y(s,t)>u }-\pk{\sup_{(s,t)\in (F_k\cup F_{k+1})\times [-\delta, \delta]}Y(s,t)>u },
\EQN
and using Theorem \ref{Th1} we have, as $u\rw\IF$,
\BQN\label{nnneq}
\Sigma_6(u)&\leq&\Theta(u)\left(\sum_{0\leq k\leq n-1}a_1(\epsilon)\left(\frac{s_{k+1}-s_k}{|g(\la{f}(s_k))f'(\la{f}(s_k))|}+\frac{s_{k+2}-s_{k+1}}{|g(\la{f}
(s_{k+1}))f'(\la{f}(s_{k+1}))|}\right)\right.\nonumber\\
&&\quad \left.-\sum_{0\leq k\leq n-1}a_1(-\epsilon)\frac{s_{k+2}-s_k}{|g(\la{f}(s_k))f'(\la{f}(s_k))|}\right)\nonumber\\
&:=&b(n,\epsilon)\Theta(u).
\EQN

Observe that
\BQNY
\pk{\sup_{(s,t)\in F_k\times [-\delta, \delta]}Y(s,t)>u, \sup_{(s,t)\in F_l\times [-\delta, \delta]}Y(s,t)>u }\leq \pk{\sup_{(s,t,s',t')\in F_k\times [-\delta, \delta]\times F_l\times [-\delta, \delta]}Y(s,t)+Y(s',t')>2u }.
\EQNY
Since for any $n$ and  $\delta$, there exists $0<\epsilon(n,\delta)<1$ such that
$$\sup_{0\leq k+1<l\leq n}\sup_{(s,t,s',t')\in F_k\times [-\delta, \delta]\times F_l\times [-\delta, \delta]} Var\left(Y(s,t)+Y(s',t')\right)<4-\epsilon(n,\delta),$$
then  Borell-TIS inequality (\cite{AdlerTaylor}) leads to
\BQN\label{eq4}
 \Sigma_7(u)\leq \sum_{0\leq k+1<l\leq n} e^{-\frac{(2u-2a)^2}{2(4-\epsilon(n,\delta))}}\leq n^2e^{-\frac{(2u-2a)^2}{2(4-\epsilon(n,\delta))}}=o\left(\Theta(u)\right), \quad u\rw\IF,
\EQN
 with $a=\E{\sup_{(s,t)\in E_1}Y(s,t)}$. Note that
\BQNY
\lim_{\epsilon\rw 0}\lim_{n\rw\IF} b(n,\epsilon)=\lim_{\epsilon\rw 0}2\left(a_1(\epsilon)-a_1(-\epsilon)\right)
\int_{f(T_1)+\delta}^{f(T_2)-\delta}\left|g(\la{f}(s))f'(\la{f}(s))\right|^{-1}ds=0.
\EQNY
Therefore, we  conclude that
\BQN\label{asym1}
\int_{f(T_1)+\delta}^{f(T_2)-\delta}\left|g(\la{f}(s))f'(\la{f}(s))\right|^{-1}ds&\leq&\lim_{\epsilon\rw 0}\lim_{n\rw\IF}\liminf_{u\rw\IF}
\frac{\pi_{4,\delta}(u)}{\Theta(u)}\nonumber\\
&\leq&\lim_{\epsilon\rw 0}\lim_{n\rw\IF}\limsup_{u\rw\IF}\frac{\pi_{4,\delta}(u)}{\Theta(u)}\nonumber\\
&&\leq \int_{f(T_1)+\delta}^{f(T_2)-\delta}\left|g(\la{f}(s))f'(\la{f}(s))\right|^{-1}ds
\EQN
Next we focus on  $\pi_{5,\delta}(u)$.
By (\ref{cor4}) and mean value theorem we have, for $\delta$ small enough,
\BQNY
1-r_Y(s,t,s',t')\leq 2\rho_1^2(|s-s'|)+2\rho_2^2\left(|\la{f}(s-t)-\la{f}(s',t')|\right)\leq 2\rho_1^2(|s-s'|)+2\rho_2^2\left(\frac{|s-s'-(t-t')|}{|f'(\la{f}(\theta))|}\right),
\EQNY
with $(s,t), (s',t')\in E_1\cap E_3$ and $\theta\in (s'-t', s-t)$. Note that {\bf F} implies that there exists $C>0$ such that $\frac{1}{|f'(\la{f}(\theta))|}\leq C$, which leads to, for $\delta$ small enough,
\BQNY
1-r_Y(s,t,s',t')\leq 2\rho_1^2(|s-s'|)+3C^{\alpha_2}\rho_2^2\left(|s-s'-(t-t')|\right), \quad (s,t), (s',t')\in E_1\cap E_3.
\EQNY
Further, by  Lemma 6.4 in \cite{KEP20162}, we have for $\delta$ small enough,
\BQN\label{nnneq1}
1-r_Y(s,t,s',t')&\leq& 4\rho_1^2(|s-s'|)+4C^{\alpha_2}\rho_2^2\left(|t-t'|\right)\nonumber\\
&\leq& 1-r_Z(8^{1/\alpha_1}s,8^{1/\alpha_2}Ct,8^{1/\alpha_1}s',8^{1/\alpha_2}Ct'), \quad (s,t), (s',t')\in E_1\cap E_3.
\EQN
Moreover, by (\ref{var4}), for $\delta$ sufficiently small
\BQN\label{nnneq2}
1-\sigma_Y(s,t)\geq \frac{c_1^\beta}{2} v^2(|t|), \quad s,t\in E_1\cap E_3.
\EQN
Thus, by Slepian inequality, for $\delta$ small enough,
\BQN\label{eq5}
\pi_{5,\delta}(u)&\leq& \pk{\sup_{(s,t)\in E_1\cap E_{3}}\frac{Z(8^{1/\alpha_1}s,8^{1/\alpha_2}Ct)}{1+\frac{c_1^\beta}{4} v^2(|t|)}>u}\nonumber\\
&\leq& \pk{\sup_{(s,t)\in  E_{3}}\frac{Z(8^{1/\alpha_1}s,8^{1/\alpha_2}Ct)}{1+\frac{c_1^\beta}{4} v^2(|t|)}>u}\nonumber\\
&=&\pk{\sup_{(s,t)\in [-\delta,\delta]^2}\frac{Z(8^{1/\alpha_1}s,8^{1/\alpha_2}Ct)}{1+\frac{c_1^\beta}{4} v^2(|t|)}>u}.
\EQN
Further, applying Theorem \ref{Th1}, we have
\BQN\label{asym2}
\pi_{5,\delta}(u)\leq 2\delta 8^{1/\alpha_1+1/\alpha_2} 4^{1/\beta}c_1^{-1} C\Theta(u)(1+o(1)), \quad u\rw\IF.
\EQN
Similarly,
\BQN\label{asym3}
\pi_{6,\delta}(u)\leq 2\delta 8^{1/\alpha_1+1/\alpha_2} 4^{1/\beta}c_1^{-1} C\Theta(u)(1+o(1)), \quad u\rw\IF.
\EQN
Finally, we focus on $\pk{\sup_{(s,t)\in E_1\setminus \{(s,t):|t|<\delta\}}Y(s,t)>u}$.
By  Borell-TIS inequality (\cite{AdlerTaylor}) and the fact that
$$\sup_{(s,t)\in  E_1\setminus \{(s,t): |t|<\delta\} }\sigma_Y(s,t)\leq 1-\epsilon, \quad \hbox{with }\epsilon>0, $$
 we have
\BQN\label{asym4}
\pk{\sup_{(s,t)\in  E_1\setminus \{(s,t): |t|<\delta\}}Y(s,t)>u}\leq 2e^{-\frac{(u-a)^2}{2(1-\epsilon)^2}}=o\left(\Theta(u)\right), \quad u\rw\IF,
\EQN
with $a=\E{\sup_{(s,t)\in E_1}Y(s,t)}$.
Inserting (\ref{asym1}), (\ref{asym2})-(\ref{asym4}) into (\ref{main4}), letting  $u \rw\IF$, $n\rw\IF$, $\delta\rw 0$ and $\epsilon\rw 0$ in turn, we have
\BQNY
\pk{\sup_{(s,t)\in E_1}Y(s,t)>u}&\sim& \int_{f(T_1)}^{f(T_2)}\left|g(\la{f}(s))f'(\la{f}(s))\right|^{-1}ds\Theta(u)\\
&\sim & \int_{T_1}^{T_2}\left|g(t)\right|^{-1}dt\Theta(u).
\EQNY
{\it \underline{Case $\gamma_2\in (0,\IF)$}}. By  Theorem \ref{Th1}, we have that
\BQNY
\pk{\sup_{(s,t)\in F_k\times [-\delta, \delta]}\frac{Z((1+\epsilon)s,(1+\epsilon)(f'(\la{f}(s_k)))^{-1}t)}{1+(1-\epsilon)|g(\la{f}(s_k))|^\beta v^2(|t|)}>u}
 \sim (1+\epsilon)\mathcal{H}_{\alpha_1} \mathcal{P}_{\alpha_2}^{b(\epsilon, s_k)}(s_{k+1}-s_k)
\frac{\Psi(u)}{\overleftarrow{\LL}_1(1/u)},
\EQNY
 with $b(\epsilon, s)=(1-\epsilon)(1+\epsilon)^{-\alpha_2}\gamma_2|g(\la{f}(s))f'(\la{f}(s))|^{\alpha_2}$.
Further, as $u\rw\IF, n\rw\IF$,
\BQNY
\sum_{k=0}^{n}\pk{\sup_{(s,t)\in F_k\times [-\delta, \delta]}Y(s,t)>u}
 &\leq& (1+\epsilon)\mathcal{H}_{\alpha_1}\sum_{k=0}^{n}\mathcal{P}_{\alpha_2}^{b(\epsilon, s_k)}(s_{k+1}-s_k)
\frac{\Psi(u)}{\overleftarrow{\LL}_1(1/u)}(1+o(1))\\
&\leq& (1+\epsilon)\mathcal{H}_{\alpha_1}\int_{f(T_1)+\delta}^{f(T_2)-\delta}\mathcal{P}_{\alpha_2}^{b(\epsilon, s)}ds
\frac{\Psi(u)}{\overleftarrow{\LL}_1(1/u)}(1+o(1)).
\EQNY
Similarly, as $u\rw\IF, n\rw\IF$,
\BQNY
\sum_{k=0}^{n}\pk{\sup_{(s,t)\in F_k\times [-\delta, \delta]}Y(s,t)>u}
\geq (1-\epsilon)\mathcal{H}_{\alpha_1}\int_{f(T_1)+\delta}^{f(T_2)-\delta}\mathcal{P}_{\alpha_2}^{b(-\epsilon, s)}ds
\frac{\Psi(u)}{\overleftarrow{\LL}_1(1/u)}(1+o(1)).
\EQNY
By (\ref{eq3}) and Theorem \ref{Th1}, we have, as $u\rw\IF$ and $n\rw\IF$,
\BQNY
\Sigma_6(u)\leq \mathcal{H}_{\alpha_1}\frac{\Psi(u)}{\overleftarrow{\LL}_1(1/u)}\sum_{0\leq k\leq n-1}\left((1+\epsilon)\mathcal{P}_{\alpha_2}^{b(\epsilon, s_k)}(s_{k+1}-s_k)+(1+\epsilon)\mathcal{P}_{\alpha_2}^{b(\epsilon, s_{k+1})}(s_{k+2}-s_{k+1})-(1-\epsilon)\mathcal{P}_{\alpha_2}^{b(-\epsilon, s_{k})}(s_{k+2}-s_{k})\right).
\EQNY
By the continuity of $\mathcal{P}_{\alpha_2}^{\gamma}$ with respect to $\gamma$ for $\gamma\in (0,\IF)$, we have that
\BQNY
&&\lim_{\epsilon\rw 0}\lim_{n\rw\IF}\sum_{0\leq k\leq n-1}\left((1+\epsilon)\mathcal{P}_{\alpha_2}^{b(\epsilon, s_k)}(s_{k+1}-s_k)+(1+\epsilon)\mathcal{P}_{\alpha_2}^{b(\epsilon, s_{k+1})}(s_{k+2}-s_{k+1})-(1-\epsilon)\mathcal{P}_{\alpha_2}^{b(-\epsilon, s_{k})}(s_{k+2}-s_{k})\right)\\
&&\quad =2\int_{f(T_1)+\delta}^{f(T_2)-\delta}\mathcal{P}_{\alpha_2}^{b(0, s)}ds-2\int_{f(T_1)+\delta}^{f(T_2)-\delta}\mathcal{P}_{\alpha_2}^{b(0, s)}ds=0.
\EQNY
Using the same argument as given in (\ref{eq4}), we have that
\BQNY
 \Sigma_7(u)=o\left(\frac{\Psi(u)}{\overleftarrow{\LL}_1(1/u)}\right), \quad u\rw\IF.
\EQNY
Therefore, in view of (\ref{eq6}) we conclude that
\BQNY
\lim_{\epsilon\rw 0}\lim_{n\rw\IF}\lim_{u\rw\IF}\frac{\pi_{4,\delta}(u)}{\frac{\Psi(u)}{\overleftarrow{\LL}_1(1/u)}}=
\mathcal{H}_{\alpha_1}\int_{f(T_1)+\delta}^{f(T_2)-\delta}\mathcal{P}_{\alpha_2}^{\gamma_2|g(\la{f}(s))f'(\la{f}(s))|^{\alpha_2}}ds.
\EQNY
By (\ref{eq5}) and Theorem \ref{Th1}, we have that for $\delta$ sufficiently small,
\BQNY
\pi_{5,\delta}(u)&\leq& \pk{\sup_{(s,t)\in [-\delta,\delta]^2}\frac{Z(8^{1/\alpha_1}s,8^{1/\alpha_2}Ct)}{1+\frac{c_1^\beta}{4} v^2(|t|)}>u}\\
&\leq&2\delta 8^{1/\alpha_1}\mathcal{H}_{\alpha_1}\mathcal{P}_{\alpha_2}^{2\gamma_2c_1^{\alpha_2} C^{\alpha_2}}\frac{\Psi(u)}{\overleftarrow{\LL}_1(1/u)}(1+o(1)), \quad u\rw\IF.
\EQNY
Similarly, for $\delta$ small enough,
\BQNY
\pi_{6,\delta}(u)
\leq2\delta 8^{1/\alpha_1}\mathcal{H}_{\alpha_1}\mathcal{P}_{\alpha_2}^{2\gamma_2c_1^{\alpha_2} C^{\alpha_2}}\frac{\Psi(u)}{\overleftarrow{\LL}_1(1/u)}(1+o(1)), \quad u\rw\IF.
\EQNY
Thus it follows from the above asymptotics and  (\ref{main4}) that
\BQNY
\liminf_{u\rw\IF}\frac{\pk{\sup_{(s,t)\in E_1}Y(s,t)>u}}{\frac{\Psi(u)}{\overleftarrow{\LL}_1(1/u)}}&\geq& \mathcal{H}_{\alpha_1}\int_{f(T_1)+\delta}^{f(T_2)-\delta}\mathcal{P}_{\alpha_2}^{\gamma_2|g(\la{f}(s))f'(\la{f}(s))|^{\alpha_2}}ds,\\ \limsup_{u\rw\IF}\frac{\pk{\sup_{(s,t)\in E_1}Y(s,t)>u}}{\frac{\Psi(u)}{\overleftarrow{\LL}_1(1/u)}}&\leq& \mathcal{H}_{\alpha_1}\int_{f(T_1)+\delta}^{f(T_2)-\delta}\mathcal{P}_{\alpha_2}^{\gamma_2|g(\la{f}(s))f'(\la{f}(s))|^{\alpha_2}}ds
+4\delta\mathcal{H}_{\alpha_1}\mathcal{P}_{\alpha_2}^{2\gamma_2c_1^{\alpha_2} C^{\alpha_2}}.
\EQNY
Letting $\delta\rw 0$, we have
\BQNY
\pk{\sup_{(s,t)\in E_1}Y(s,t)>u}&\sim&\mathcal{H}_{\alpha_1}\frac{\Psi(u)}{\overleftarrow{\LL}_1(1/u)}\int_{f(T_1)}^{f(T_2)}
\mathcal{P}_{\alpha_2}^{\gamma_2|g(\la{f}(s))f'(\la{f}(s))|^{\alpha_2}}ds\\
&\sim&\mathcal{H}_{\alpha_1}\frac{\Psi(u)}{\overleftarrow{\LL}_1(1/u)}\int_{T_1}^{T_2}
\mathcal{P}_{\alpha_2}^{\gamma_2|g(t)f'(t)|^{\alpha_2}}|f'(t)|dt.
\EQNY
{\it \underline{Case $\gamma_2=\IF$}}.  Replacing Piterbarg constants by $1$ in the proof of case $\gamma_2\in (0,\IF)$, we can establish the claim. \QED\\
\prooftheo{th32}  Without loss of generality, we assume that $f'(t)>0, t\in(T_1, T_2)$.
 Let $Y_1(s,t)=X(s,\la{f}(t)), (s,t)\in F_1:=\{(s,t): S_1\leq s\leq S_2, f(T_1)\leq t\leq f(T_2)\}$. Then
\BQNY
\pk{\sup_{(s,t)\in [S_1,S_2]\times [T_1,T_2] }X(s,t)>u}=\pk{\sup_{(s,t)\in F_1}Y_1(s,t)>u}.
\EQNY
Note that $\sigma_{Y_1}(s,t), (s,t)\in F_1$ attains its maximum over the line $\{(t,t), f(T_1)\leq t\leq f(T_2)\}$ and  satisfies
\BQNY
1-\sigma_{Y_1}(s,t)\sim v^2(g(\la{f}(t))|s-t|), \quad |s-t|\rw 0, (s,t)\in F_1.
\EQNY
Moreover,
\BQNY
1-r_{Y_1}(s,t,s',t')\sim \rho_1^2(|s-s'|)+\rho_2^2\left(\frac{|t-t'|}{f'(\la{f}(t))}\right),\quad |s-s'|, |t-t'|, |s-t|\rw 0, (s',t'), (s,t)\in F_1.
\EQNY
Let $F_2=\{(s,t): |s-t|<\delta,  f(T_1)+\delta\leq t\leq f(T_2)-\delta\}$, $F_3=\{(s,t): |s-t|<\delta, f(T_1)\leq t\leq f(T_1)+\delta\}$, $F_4=\{(s,t): |s-t|<\delta, f(T_2)-\delta\leq t\leq f(T_2)\}$. We have
\BQN\label{main10}
\pi_{7,\delta}(u)\leq \pk{\sup_{(s,t)\in F_1}Y_1(s,t)>u}\leq \sum_{i=7}^9\pi_{i,\delta}(u)+\pk{\sup_{(s,t)\in F_1\setminus \{(s,t):|s-t|\leq\delta\}}Y_1(s,t)>u},
\EQN
where $\pi_{i,\delta}(u)=\pk{\sup_{(s,t)\in F_{i-5}\cap F_1}Y_1(s,t)>u}, i=7,8,9.$
Let $G_k=\{(s,t): |s-t|<\delta,  t_k\leq t\leq t_{k+1}\}$ with $t_k=f(T_1)+\delta+\frac{k(f(T_2)-f(T_1)-2\delta)}{n}, 0\leq k\leq n$.
We have
\BQN\label{main11}
\sum_{k=0}^{n}\pk{\sup_{(s,t)\in G_k }Y_1(s,t)>u}-\sum_{i=8}^9\Sigma_i(u)\leq\pi_{7,\delta}(u)\leq \sum_{k=0}^{n}\pk{\sup_{(s,t)\in G_k }Y_1(s,t)>u},
\EQN
where
\BQNY
\Sigma_8(u)&=&\sum_{0\leq k\leq n-1}\pk{\sup_{(s,t)\in G_k }Y_1(s,t)>u, \sup_{(s,t)\in G_{k+1} }Y_1(s,t)>u},\\
\Sigma_9(u)&=&\sum_{0\leq k<l\leq n, l\geq k+2}\pk{\sup_{(s,t)\in G_k }Y_1(s,t)>u, \sup_{(s,t)\in G_{l} }Y_1(s,t)>u}.
\EQNY
Analogously as in the proof of Theorem \ref{th31}, let $Z(s,t)$ be a homogeneous Gaussian field with variance $1$ and correlation function satisfying (\ref{corr}). Then it follows that for $0<\epsilon<1$ sufficiently small, if $\delta$ sufficiently small and $n$ sufficiently large, then
\BQNY
1-r_{Y_1}(s,t,s',t')&\geq& 1-r_Z\left((1-\epsilon)s,(1-\epsilon)(f'(\la{f}(t_k)))^{-1}t,(1-\epsilon)s',(1-\epsilon)(f'(\la{f}(t_k)))^{-1}t'\right),\\
1-r_{Y_1}(s,t,s',t')&\leq& 1-r_Z\left((1+\epsilon)s,(1+\epsilon)(f'(\la{f}(t_k)))^{-1}t,(1+\epsilon)s',(1+\epsilon)(f'(\la{f}(t_k)))^{-1}t'\right)
\EQNY
hold for $(s,t),(s',t')\in G_k, 0\leq k\leq n$, and
\BQNY
(1-\epsilon) (g(\la{f}(t_k)))^\beta v^2(|s-t|)\leq 1-\sigma_{Y_1}(s,t)\leq (1+\epsilon) (g(\la{f}(t_k)))^\beta v^2(|s-t|)
\EQNY
holds for $(s,t)\in G_k, 0\leq k\leq n $.
 By Slepian inequality, we have
\BQNY
\pk{\sup_{(s,t)\in G_k }\frac{Z((1-\epsilon)s,(1-\epsilon)(f'(\la{f}(t_k)))^{-1}t)}{1+(1+\epsilon) (g(\la{f}(t_k)))^\beta v^2(|s-t|)}>u}&\leq&\pk{\sup_{(s,t)\in G_k }Y_1(s,t)>u}\\&\leq& \pk{\sup_{(s,t)\in G_k }\frac{Z((1+\epsilon)s,(1+\epsilon)(f'(\la{f}(t_k)))^{-1}t)}{1+(1-\epsilon) (g(\la{f}(t_k)))^\beta v^2(|s-t|)}>u}
\EQNY
Direct calculation shows that
\BQNY
&&1-r_Z\left((1\pm\epsilon)s,(1\pm\epsilon)(f'(\la{f}(t_k)))^{-1}t,(1\pm\epsilon)s',(1\pm\epsilon)(f'(\la{f}(t_k)))^{-1}t'\right)\\
&&\quad\sim (1\pm\epsilon)^{\alpha_1}\rho_1^2(|s-s'|)+(1\pm\epsilon)^{\alpha_2}(f'(\la{f}(t_k)))^{-\alpha_2}\rho_2^2(|t-t'|), \quad |s-s'|, |t-t'|\rw 0, |s-t|\rw 0,
\EQNY
for $(s,t), (s',t')\in G_k, 0\leq k\leq n$.\\
\underline{\it Case $\gamma_1=0$}. By Theorem \ref{Th2}, we have
\BQN\label{asy10}
\pk{\sup_{(s,t)\in G_k }\frac{Z((1\pm\epsilon)s,(1\pm\epsilon)(f'(\la{f}(t_k)))^{-1}t)}{1+(1\mp\epsilon) (g(\la{f}(t_k)))^\beta v^2(|s-t|)}>u}\sim a_1(\pm\epsilon)\frac{t_{k+1}-t_k}{|g(\la{f}(t_k))f'(\la{f}(t_k))|}\Theta(u), \quad u\rw\IF,
\EQN
where $\Theta(u)$ is defined in the proof of Theorem \ref{th31} and $a_1(\pm \epsilon)=(1\mp\epsilon)^{-1/\beta}(1\pm\epsilon)^2$. This implies that, as $u\rw\IF$,
\BQN\label{asy3}
a_1(-\epsilon)\sum_{k=0}^{n}\frac{t_{k+1}-t_k}{|g(\la{f}(t_k))f'(\la{f}(t_k))|}\leq \frac{\sum_{k=0}^{n}\pk{\sup_{(s,t)\in G_k }Y_1(s,t)>u}}{\Theta(u)}\leq a_1(\epsilon)\sum_{k=0}^{n}\frac{t_{k+1}-t_k}{|g(\la{f}(t_k))f'(\la{f}(t_k))|}
\EQN
Using the same arguments as given in  (\ref{nnneq}) and (\ref{eq4}), we have that, as $u\rw\IF$,
\BQN
\Sigma_8(u)\leq b_1(n,\epsilon)\Theta(u), \quad \Sigma_9(u)=o\left(\Theta(u)\right),
\EQN
where $ b_1(n,\epsilon)>0$ is function of $n$ and $\epsilon$ such that $\lim_{\epsilon\rw 0}\lim_{n\rw\IF}b_1(n,\epsilon)=0$.
Similarly as in (\ref{nnneq1})- (\ref{eq5}), we have
for $i=8,9$,
\BQNY
\pi_{i,\delta}(u)&\leq& \pk{\sup_{(s,t)\in F_{i-5}\cap F_{1}}\frac{Z(8^{1/\alpha_1}s,8^{1/\alpha_2}Ct)}{1+\frac{c_1^\beta}{4} v^2(|s-t|)}>u}\nonumber\\
&\leq& \pk{\sup_{(s,t)\in   F_{i-5}}\frac{Z(8^{1/\alpha_1}s,8^{1/\alpha_2}Ct)}{1+\frac{c_1^\beta}{4} v^2(|s-t|)}>u}\nonumber\\
&=&\pk{\sup_{(s,t)\in [-2\delta,2\delta]^2}\frac{Z(8^{1/\alpha_1}s,8^{1/\alpha_2}Ct)}{1+\frac{c_1^\beta}{4} v^2(|s-t|)}>u}.
\EQNY
Further, applying Theorem \ref{Th2}, we have
\BQN\label{asy4}
\pi_{8,\delta}(u)&\leq& 4\delta 8^{2/\alpha} 4^{1/\beta}c_1^{-1} C\Theta(u)(1+o(1)),\nonumber\\
\pi_{9,\delta}(u)&\leq& 4\delta 8^{2/\alpha} 4^{1/\beta}c_1^{-1} C\Theta(u)(1+o(1)), \quad u\rw\IF.
\EQN
 Combination of (\ref{main10})-(\ref{asy4}) leads to
 \BQNY
\limsup_{u\rw\IF}\frac{\pk{\sup_{(s,t)\in F_1}Y_1(s,t)>u}}{\Theta(u)}&\leq& a_1(\epsilon)\sum_{k=0}^{n}\frac{t_{k+1}-t_k}{|g(\la{f}(t_k))f'(\la{f}(t_k))|}+8\delta 8^{2/\alpha} 4^{1/\beta}c_1^{-1}C,\\
\liminf_{u\rw\IF}\frac{\pk{\sup_{(s,t)\in F_1}Y_1(s,t)>u}}{\Theta(u)}&\geq& a_1(-\epsilon)\sum_{k=0}^{n}\frac{t_{k+1}-t_k}{|g(\la{f}(t_k))f'(\la{f}(t_k))|}- b_1(n,\epsilon).
 \EQNY
 Letting $n\rw\IF$, $\epsilon\rw 0$ and $\delta\rw 0$ in turn, we have
 \BQNY
 \pk{\sup_{(s,t)\in F_1}Y_1(s,t)>u}&\sim& \int_{f(T_1)}^{f(T_2)}\left|g(\la{f}(s))f'(\la{f}(s))\right|^{-1}ds\Theta(u)\\
&\sim & \int_{T_1}^{T_2}\left|g(t)\right|^{-1}dt\Theta(u).
 \EQNY
 \underline{\it Case $\gamma_1\in (0,\IF)$}.
  By Theorem \ref{Th2}, we have
\BQNY
\pk{\sup_{(s,t)\in G_k }\frac{Z((1\pm\epsilon)s,(1\pm\epsilon)(f'(\la{f}(t_k)))^{-1}t)}{1+(1\mp\epsilon) (g(\la{f}(t_k)))^\beta v^2(|s-t|)}>u}&\sim& \mathcal{H}_{\alpha_1}^{(1\mp\epsilon)(1\pm\epsilon)^{-\alpha_1}\gamma_1(g(\la{f}(t_k)))^{\alpha_1},-\eta^{-1/\alpha_1}|f'(\la{f}(t_k))|}\\
&&  \times\frac{(1\pm\epsilon)(t_{k+1}-t_k)}{|f'(\la{f}(t_k))|}\frac{\Psi(u)}{\la{\rho}_2(1/u)}, \quad u\rw\IF,
\EQNY
 for $0\leq k\leq n$. Hence, as $u\rw\IF$,
 \BQN\label{eq10}
 \frac{\sum_{k=0}^{n}\pk{\sup_{(s,t)\in G_k }Y_1(s,t)>u}}{\frac{\Psi(u)}{\la{\rho}_2(1/u)}}&\leq& \sum_{k=0}^{n}\mathcal{H}_{\alpha_1}^{(1-\epsilon)(1+\epsilon)^{-\alpha_1}\gamma_1(g(\la{f}(t_k)))^{\alpha_1},-\eta^{-1/\alpha_1}|f'(\la{f}(t_k))|}
\frac{(1+\epsilon)(t_{k+1}-t_k)}{|f'(\la{f}(t_k))|},\nonumber\\
\frac{\sum_{k=0}^{n}\pk{\sup_{(s,t)\in G_k }Y_1(s,t)>u}}{\frac{\Psi(u)}{\la{\rho}_2(1/u)}}&\geq& \sum_{k=0}^{n
}\mathcal{H}_{\alpha_1}^{(1+\epsilon)(1-\epsilon)^{-\alpha_1}\gamma_1(g(\la{f}(t_k)))^{\alpha_1},-\eta^{-1/\alpha_1}|f'(\la{f}(t_k))|}
\frac{(1-\epsilon)(t_{k+1}-t_k)}{|f'(\la{f}(t_k))|}.
\EQN
 Using the same arguments as in (\ref{eq3}), by Theorem \ref{Th2} we have that
 \BQN
 \Sigma_8(u)\leq b_2(n,\epsilon) \frac{\Psi(u)}{\la{\rho}_2(1/u)}, \quad u\rw\IF,
 \EQN
 with
 \BQNY
 b_2(n,\epsilon)&=&2(1+\epsilon)\sum_{k=0}^{n-2}\mathcal{H}_{\alpha_1}^{(1-\epsilon)(1+\epsilon)^{-\alpha_1}\gamma_1(g(\la{f}(t_k)))^{\alpha_1},-\eta^{-1/\alpha_1}|f'(\la{f}(t_k))|}
\frac{(t_{k+1}-t_k)}{|f'(\la{f}(t_k))|}\\
 &&-(1-\epsilon)\sum_{k=0}^{n-2}\mathcal{H}_{\alpha_1}^{(1+\epsilon)(1-\epsilon)^{-\alpha_1}\gamma_1(g(\la{f}(t_k)))^{\alpha_1},-\eta^{-1/\alpha_1}
 |f'(\la{f}(t_k))|}
\frac{(t_{k+2}-t_k)}{|f'(\la{f}(t_k))|}.
 \EQNY
 Note that
 \BQNY
\lim_{\epsilon\rw 0} \lim_{n\rw\IF}b_2(n,\epsilon)&=&\lim_{\epsilon\rw 0} 2(1+\epsilon)\int_{f(T_1)+\delta}^{f(T_2)-\delta}\mathcal{H}_{\alpha_1}^{(1-\epsilon)(1+\epsilon)^{-\alpha_1}\gamma_1(g(\la{f}(t)))^{\alpha_1},-\eta^{-1/\alpha_1}|f'(\la{f}(t))|}
\frac{1}{|f'(\la{f}(t))|}dt\\
&& -\lim_{\epsilon\rw 0}2(1-\epsilon)\int_{f(T_1)+\delta}^{f(T_2)-\delta}\mathcal{H}_{\alpha_1}^{(1+\epsilon)(1-\epsilon)^{-\alpha_1}\gamma_1(g(\la{f}(t)))^{\alpha_1},-\eta^{-1/\alpha_1}|f'(\la{f}(t))|}
\frac{1}{|f'(\la{f}(t))|}dt=0.
 \EQNY
 Similarly as in  (\ref{eq4}), we have that, as $u\rw\IF$,
\BQN\label{eq11}
 \Sigma_9(u)=o\left( \frac{\Psi(u)}{\la{\rho}_2(1/u)}\right).
\EQN
Analogously as in (\ref{asy4}), we derive that
 \BQNY
\pi_{8,\delta}(u)&\leq& \mathbb{Q}\delta \frac{\Psi(u)}{\la{\rho}_2(1/u)}(1+o(1)),\nonumber\\
\pi_{9,\delta}(u)&\leq& \mathbb{Q}\delta  \frac{\Psi(u)}{\la{\rho}_2(1/u)}(1+o(1)), \quad u\rw\IF,
\EQNY
which combined with (\ref{main10})-(\ref{main11}) and (\ref{eq10})-(\ref{eq11}) lead to
\BQNY
\limsup_{u\rw\IF}\frac{\pk{\sup_{(s,t)\in F_1}Y_1(s,t)>u}}{\frac{\Psi(u)}{\la{\rho}_2(1/u)}}&\leq& \sum_{k=0}^{n
}\mathcal{H}_{\alpha_1}^{(1-\epsilon)(1+\epsilon)^{-\alpha_1}\gamma_1(g(\la{f}(t_k)))^{\alpha_1},-\eta^{-1/\alpha_1}|f'(\la{f}(t_k))|}
\frac{(1+\epsilon)(t_{k+1}-t_k)}{|f'(\la{f}(t_k))|}+2\mathbb{Q}\delta,\\
\liminf_{u\rw\IF}\frac{\pk{\sup_{(s,t)\in F_1}Y_1(s,t)>u}}{\frac{\Psi(u)}{\la{\rho}_2(1/u)}}&\geq& \sum_{k=0}^{n}\mathcal{H}_{\alpha_1}^{(1+\epsilon)(1-\epsilon)^{-\alpha_1}\gamma_1(g(\la{f}(t_k)))^{\alpha_1},-\eta^{-1/\alpha_1}|f'(\la{f}(t_k))|}
\frac{(1-\epsilon)(t_{k+1}-t_k)}{|f'(\la{f}(t_k))|}-b_2(n,\epsilon).
\EQNY
 Letting $n\rw\IF$, $\delta\rw 0$ and $\epsilon\rw 0$ respectively in the above inequalities, we derive that
\BQNY
\pk{\sup_{(s,t)\in F_1}Y_1(s,t)>u}&\sim& \int_{f(T_1)}^{f(T_2)}\mathcal{H}_{\alpha_1}^{\gamma_1(g(\la{f}(t)))^{\alpha_1},-\eta^{-1/\alpha_1}|f'(\la{f}(t))|}
\frac{1}{|f'(\la{f}(t))|}dt \frac{\Psi(u)}{\la{\rho}_2(1/u)}\\
&\sim& \int_{T_1}^{T_2}\mathcal{H}_{\alpha_1}^{\gamma_1(g(t))^{\alpha_1},-\eta^{1/\alpha_1}|f'(t)|}
dt \frac{\Psi(u)}{\la{\rho}_2(1/u)}, \quad u\rw\IF.
\EQNY
 \underline{\it Case $\gamma_1=\IF$}. By Theorem \ref{Th2},
 replacing $\mathcal{H}_{\alpha_1}^{x,y}$ by $(|y|^{\alpha_1}+1)^{1/\alpha_1}\mathcal{H}_{\alpha_1}$ in the proof of case $\gamma_1\in(0,\IF)$, we have that
 \BQNY
 \pk{\sup_{(s,t)\in F_1}Y_1(s,t)>u}&\sim& \int_{T_1}^{T_2}(\eta^{-1}|f'(t)|^{\alpha_1}+1)^{1/\alpha_1}dt\mathcal{H}_{\alpha_1} \frac{\Psi(u)}{\la{\rho}_2(1/u)}, \quad u\rw\IF.
 \EQNY
 This completes the proof. \QED\\
 \prooftheo{th33} Without loss of generality, we assume that $f'(t)>0, t\in (T_1, T_2)$.
 Let $Y_2(s,t)=X(s+t,\la{f}(t)), (s,t)\in F_1:=\{(s,t): S_1\leq s+t\leq S_2, f(T_1)\leq t\leq f(T_2)\}$. Then
\BQNY
\pk{\sup_{(s,t)\in [S_1,S_2]\times [T_1,T_2] }X(s,t)>u}=\pk{\sup_{(s,t)\in F}Y_2(s,t)>u}.
\EQNY
Clearly, $\sigma_Y(s,t), (s,t)\in F_1$ attains its maximum on the line $\{(0,t): f(T_1)\leq t\leq f(T_2)\}$. Moreover, by (\ref{cor4}), (\ref{var4}) and Lemma 6.4 in \cite{KEP20162}, we have for $(s,t), (s',t')\in F_1$
\BQNY
1-r_Y(s,t,s',t')\sim \rho_1^2(|s-s'+t-t'|)+\rho_2^2\left(\frac{|t-t'|}{f'(\la{f}(t))}\right)\sim  \rho_1^2(|s-s'|)+\rho_2^2\left(\frac{|t-t'|}{f'(\la{f}(t))}\right), \quad s,s'\rw 0, |t-t'|\rw 0,
\EQNY
and
\BQNY
1-\sigma_Y(s,t)\sim v^2(g(\la{f}(t))|s|), \quad |s|\rw 0.
\EQNY
 Let $F_5=\{(s,t): |s|\leq \delta, f(T_1)+\delta\leq t\leq f(T_2)-\delta\}$, $F_6=\{(s,t): |s|\leq \delta, f(T_1)-\delta\leq t\leq f(T_1)+\delta\}$,
 and $F_7=\{(s,t): |s|\leq \delta, f(T_2)-\delta\leq t\leq f(T_2)+\delta\}$. Observe that
\BQNY
\pi_{10,\delta}(u)\leq \pk{\sup_{(s,t)\in F_1}Y_2(s,t)>u}\leq \sum_{i=10}^{12}\pi_{i,\delta}(u)+\pk{\sup_{(s,t)\in F\setminus \{(s,t):|s|\leq\delta\}}Y_2(s,t)>u},
\EQNY
where $\pi_{i,\delta}(u)=\pk{\sup_{(s,t)\in F_{i-5}\cap F}Y_2(s,t)>u}, i=10,11,12.$ Let $G_k=\{(s,t): |s|<\delta,  t_k\leq t\leq t_{k+1}\}$ with $t_k=f(T_1)+\delta+\frac{k(f(T_2)-f(T_1)-2\delta)}{n}, 0\leq k\leq n$.
  Then
  \BQNY
  \sum_{k=0}^{n
  }\pk{\sup_{(s,t)\in G_k}Y_2(s,t)>u}-\sum_{i=10}^{11}\Sigma_{i}(u)\leq\pi_{10,\delta}(u)\leq \sum_{k=0}^{n
  }\pk{\sup_{(s,t)\in G_k}Y_2(s,t)>u},
  \EQNY
  where
\BQNY
\Sigma_{10}(u)&=&\sum_{0\leq k\leq n-1}\pk{\sup_{(s,t)\in G_k }Y_1(s,t)>u, \sup_{(s,t)\in G_{k+1} }Y_1(s,t)>u},\\
\Sigma_{11}(u)&=&\sum_{0\leq k<l\leq n-1, l\geq k+2}\pk{\sup_{(s,t)\in G_k }Y_1(s,t)>u, \sup_{(s,t)\in G_{l} }Y_1(s,t)>u}.
\EQNY
 It follows that for $0<\epsilon<1$ sufficiently small, if $\delta$ sufficiently small and $n$ sufficiently large, then
\BQNY
1-r_{Y}(s,t,s',t')&\geq& 1-r_Z\left((1-\epsilon)s,(1-\epsilon)(f'(\la{f}(t_k)))^{-1}t,(1-\epsilon)s',(1-\epsilon)(f'(\la{f}(t_k)))^{-1}t'\right)\\
1-r_{Y}(s,t,s',t')&\leq& 1-r_Z\left((1+\epsilon)s,(1+\epsilon)(f'(\la{f}(t_k)))^{-1}t,(1+\epsilon)s',(1+\epsilon)(f'(\la{f}(t_k)))^{-1}t'\right)
\EQNY
hold for $(s,t),(s',t')\in G_k, 0\leq k\leq n$, and
\BQNY
(1-\epsilon) (g(\la{f}(t_k)))^\beta v^2(|s|)\leq 1-\sigma_{Y}(s,t)\leq (1+\epsilon) (g(\la{f}(t_k)))^\beta v^2(|s|)
\EQNY
holds for $(s,t)\in G_k, 0\leq k\leq n
 $.
 By Slepian inequality, we have
\BQNY
\pk{\sup_{(s,t)\in G_k }\frac{Z((1-\epsilon)s,(1-\epsilon)(f'(\la{f}(t_k)))^{-1}t)}{1+(1+\epsilon) (g(\la{f}(t_k)))^\beta v^2(|s|)}>u}&\leq&\pk{\sup_{(s,t)\in G_k }Y_2(s,t)>u}\\&\leq& \pk{\sup_{(s,t)\in G_k }\frac{Z((1+\epsilon)s,(1+\epsilon)(f'(\la{f}(t_k)))^{-1}t)}{1+(1-\epsilon) (g(\la{f}(t_k)))^\beta v^2(|s|)}>u}
\EQNY
Direct calculation shows that
\BQNY
&&1-r_Z\left((1\pm\epsilon)s,(1\pm\epsilon)(f'(\la{f}(t_k)))^{-1}t,(1\pm\epsilon)s',(1\pm\epsilon)(f'(\la{f}(t_k)))^{-1}t'\right)\\
&&\quad\sim (1\pm\epsilon)^{\alpha_1}\rho_1^2(|s-s'|)+(1\pm\epsilon)^{\alpha_2}(f'(\la{f}(t_k)))^{-\alpha_2}\rho_2^2(|t-t'|), \quad |s-s'|, |t-t'|\rw 0,
\EQNY
for $(s,t), (s',t')\in G_k$. \\
\underline{\it Case $\gamma_1=0$}. Using the same arguments as given in (\ref{asy10})-(\ref{asy4}), we derive that
\BQNY
 \pk{\sup_{(s,t)\in F_1}Y_2(s,t)>u}\sim
\int_{T_1}^{T_2}\left|g(t)\right|^{-1}dt\Theta(u), \quad u\rw\IF.
\EQNY
 \underline{\it Case $\gamma_1\in (0,\IF)$}. By Theorem \ref{Th1}, we have that
 \BQNY
 \pk{\sup_{(s,t)\in G_k }\frac{Z((1\pm\epsilon)s,(1\pm\epsilon)(f'(\la{f}(t_k)))^{-1}t)}{1+(1\mp\epsilon) (g(\la{f}(t_k)))^\beta v^2(|s|)}>u}\sim
 \mathcal{H}_{\alpha_2} \mathcal{P}_{\alpha_1}^{b_3(\pm\epsilon, t_k)}(1\pm\epsilon)\frac{t_{k+1}-t_k}{|f'(\la{f}(t_k))|} \frac{\Psi(u)}{\la{\rho}_2(1/u)}
 \EQNY
 with $b_3(\pm\epsilon, t_k)=(1\mp\epsilon)(1\pm\epsilon)^{-\alpha_1}\gamma_1(g(\la{f}(t_k)))^\beta$.
 Hence,
 \BQN\label{asy11}
\limsup_{u\rw\IF}\frac{ \sum_{k=0}^{n
}\pk{\sup_{(s,t)\in G_k}Y_2(s,t)>u}}{\frac{\Psi(u)}
 {\la{\rho}_2(1/u)}}&\leq& (1+\epsilon)\mathcal{H}_{\alpha_2} \sum_{k=0}^{n } \mathcal{P}_{\alpha_1}^{b_3(+\epsilon, t_k)}\frac{t_{k+1}-t_k}{|f'(\la{f}(t_k))|},\nonumber\\
  \liminf_{u\rw\IF}\frac{ \sum_{k=0}^{n}\pk{\sup_{(s,t)\in G_k}Y_2(s,t)>u}}{\frac{\Psi(u)}{\la{\rho}_2(1/u)}}&\geq&(1-\epsilon)\mathcal{H}_{\alpha_2} \sum_{k=0}^{n} \mathcal{P}_{\alpha_1}^{b_3(-\epsilon, t_k)}\frac{t_{k+1}-t_k}{|f'(\la{f}(t_k))|} .
 \EQN
 Following the same arguments as in the proof of Theorem \ref{th31}, we can get that
 \BQNY
\lim_{\epsilon\rw 0}\lim_{\delta\rw 0}\lim_{n\rw\IF} \limsup_{u\rw\IF}\frac{\sum_{i=10}^{11}\Sigma_{i}(u)+\sum_{i=11}^{12}\pi_{i,\delta}(u)+\pk{\sup_{(s,t)\in F_1\setminus \{(s,t):|s|\leq\delta\}}Y_2(s,t)>u}}{\frac{\Psi(u)}
 {\la{\rho}_2(1/u)}}=0.
 \EQNY
 Thus letting $u\rw\IF$, $n\rw\IF$, $\delta\rw 0$ and $\epsilon\rw 0$ in turn, we derive that
 \BQNY
 \pk{\sup_{(s,t)\in F}Y_2(s,t)>u}&\sim& \mathcal{H}_{\alpha_2}\int_{f(T_1)}^{f(T_2)}\mathcal{P}_{\alpha_1}^{\gamma_1(g(\la{f}(t)))^\beta}\frac{1}{|f'(\la{f}(t))|}dt\frac{\Psi(u)}
 {\la{\rho}_2(1/u)}\\
 &\sim& \mathcal{H}_{\alpha_2}\int_{T_1}^{T_2}\mathcal{P}_{\alpha_1}^{\gamma_1(g(t))^\beta}dt\frac{\Psi(u)}
 {\la{\rho}_2(1/u)}.
 \EQNY
 \underline{\it Case $\gamma_1=0$}. Letting $\mathcal{P}_{\alpha}^y=1$ in the proof of case $\gamma_1\in (0,\IF)$ establishes the claim.
 This completes the proof. \QED\\

\proofprop{cor11} Note that $\sigma(s,t)=\sqrt{\Var(B_{\alpha_1}(s)+B_{\alpha_2}(t))}$ attains $1$ at
 $\mathcal{L}=\{(s,t): |s|^{\alpha_1}+|t|^{\alpha_2}=1\}$. We first analyze local behavior of variance function as $|s|^{\alpha_1}+|t|^{\alpha_2}\rw 1$.
Observe that
\BQN\label{var5}
1-\sigma(s,t)&=&1-\sqrt{|s|^{\alpha_1}+|t|^{\alpha_2}}\nonumber\\
 &\sim& \frac{1}{2}(1-|s|^{\alpha_1}-|t|^{\alpha_2})\nonumber\\
 &=&\frac{1}{2}(1-|s|^{\alpha_1})(1-|t(1-|s|^{\alpha_1})^{-1/\alpha_2}|^{\alpha_2})\nonumber\\
&\sim& \frac{\alpha_2}{2}(1-|s|^{\alpha_1})(1-|t|(1-|s|^{\alpha_1})^{-1/\alpha_2})\nonumber\\
&=&\frac{\alpha_2}{2}(1-|s|^{\alpha_1})^{1-1/\alpha_2}\left||t|-(1-|s|^{\alpha_1})^{1/\alpha_2}\right|,
\EQN
holds as $|s|^{\alpha_1}+|t|^{\alpha_2}\uparrow 1$ with $|s|\leq 1-\delta $ and $0<\delta<1$.
Similarly,
\BQN\label{var6}
1-\sigma(s,t) \sim \frac{\alpha_1}{2}(1-|t|^{\alpha_2})^{1-1/\alpha_1}\left||s|-(1-|t|^{\alpha_2})^{1/\alpha_1}\right|
\EQN
holds as $|s|^{\alpha_1}+|t|^{\alpha_2}\uparrow 1$ with $|t|\leq 1-\delta$ and  $0<\delta<1$.
Next we focus on the local behavior of  correlation function as $|s|^{\alpha_1}+|t|^{\alpha_2}\rw 1$.  We have
\BQNY
1-r(s,t,s_1,t_1)&=&1-\frac{\mathbb{E}\left((B_{\alpha_1}(s)+B_{\alpha_2}(t))(B_{\alpha_1}(s_1)+
B_{\alpha_2}(t_1))\right)}{\sqrt{|s|^{\alpha_1}+|t|^{\alpha_2}}\sqrt{|s_1|^{\alpha_1}+|t_1|^{\alpha_2}}}\\
&=&\frac{Var
\left((B_{\alpha_1}(s)-B_{\alpha_1}(s_1)+B_{\alpha_2}(t)-B_{\alpha_2}(t_1))\right)-\left(\sqrt{|s|^{\alpha_1}+|t|^{\alpha_2}}
-\sqrt{|s_1|^{\alpha_1}+|t_1|^{\alpha_2}}\right)^2}
{2\sqrt{|s|^{\alpha_1}+|t|^{\alpha_2}}\sqrt{|s_1|^{\alpha_1}+|t_1|^{\alpha_2}}}\\
&\sim& \frac{|s-s_1|^{\alpha_1}+|t-t_1|^{\alpha_2}-\left(\sqrt{|s|^{\alpha_1}+|t|^{\alpha_2}}-\sqrt{|s_1|^{\alpha_1}+|t_1|^{\alpha_2}}\right)^2}
{2}\\
&\sim&\frac{|s-s_1|^{\alpha_1}+|t-t_1|^{\alpha_2}}{2}\left(1-\frac{(\sqrt{|s|^{\alpha_1}+|t|^{\alpha_2}}-\sqrt{|s_1|^{\alpha_1}+|t_1|^{\alpha_2}})^2}
{|s-s_1|^{\alpha_1}+|t-t_1|^{\alpha_2}}\right),
\EQNY
 as $|s|^{\alpha_1}+|t|^{\alpha_2}, |s_1|^{\alpha_1}+|t_1|^{\alpha_2}\uparrow 1$ and $(s,t)\neq (s_1,t_1)$.
Moreover,
\BQNY
\frac{\left(\sqrt{|s|^{\alpha_1}+|t|^{\alpha_2}}-\sqrt{|s_1|^{\alpha_1}+|t_1|^{\alpha_2}}\right)^2}
{|s-s_1|^{\alpha_1}+|t-t_1|^{\alpha_2}}&\leq& \frac{\left(|s|^{\alpha_1}+|t|^{\alpha_2}-|s_1|^{\alpha_1}-|t_1|^{\alpha_2}\right)^2}
{|s-s_1|^{\alpha_1}+|t-t_1|^{\alpha_2}}\\
&\leq & \frac{2\left(|s|^{\alpha_1}-|s_1|^{\alpha_1})^2+2(|t|^{\alpha_2}-|t_1|^{\alpha_2}\right)^2}
{|s-s_1|^{\alpha_1}+|t-t_1|^{\alpha_2}}, \quad (s,t)\neq (s_1,t_1).
\EQNY
Assume that $|s|>|s_1|$ and let $x=\frac{s_1}{s}$. Then for $0<\epsilon<\min_{i=1,2}\{\alpha_i, 2-\alpha_i\}$,
\BQNY
\frac{(|s|^{\alpha_1}-|s_1|^{\alpha_1})^2}
{|s-s_1|^{\alpha_1}}&=& \frac{(|s|^{\alpha_1}-|s_1|^{\alpha_1})^2}
{|s-s_1|^{\alpha_1+\epsilon}}|s-s_1|^{\epsilon}\\
&\leq&|s-s_1|^{\epsilon}|s|^{\alpha_1-\epsilon}\sup_{x\in (-1,1)}\frac{(1-|x|^{\alpha_1})^2}{|1-x|^{\alpha_1+\epsilon}}\leq \mathbb{Q}|s-s_1|^{\epsilon},
\EQNY
implying that
$$\frac{\left(\sqrt{|s|^{\alpha_1}+|t|^{\alpha_2}}-\sqrt{|s_1|^{\alpha_1}+|t_1|^{\alpha_2}}\right)^2}
{|s-s_1|^{\alpha_1}+|t-t_1|^{\alpha_2}}\leq \mathbb{Q}(|s-s_1|^{\epsilon}+|t-t_1|^\epsilon), \quad (s,t)\neq (s_1,t_1).$$
Hence,
\BQN\label{cor5}
\lim_{\delta\rw 0}\sup_{(s,t)\neq (s_1,t_1), |s-s_1|\leq \delta, |t-t_1|\leq \delta, 1-\delta\leq |s|^{\alpha_1}+|t|^{\alpha_2}, |s_1|^{\alpha_1}+|t_1|^{\alpha_2}\leq 1 }\left|\frac{2(1-r(s,t,s_1,t_1))}{|s-s_1|^{\alpha_1}+|t-t_1|^{\alpha_2}}-1\right|=0.
\EQN
 Let $$\Pi_{\delta_1,\delta_2,\delta_3,\delta_4}(u)=\pk{\sup_{(s,t)\in F(\delta_1,\delta_2, \delta_3, \delta_4)}(B_{\alpha_1}(s)+B_{\alpha_2}(t))>u}$$
 with
 $F(\delta_1,\delta_2, \delta_3, \delta_4)=\{(s,t): |s|^{\alpha_1}+|t|^{\alpha_2}\leq 1, \delta_1\leq s\leq \delta_2, \delta_3 \leq t\leq \delta_4\}$.
 Then
 by the symmetry of $B_{\alpha}(s)+B_{\alpha_2}(t)$, we have that for any $0<\delta<1$ and $0<\epsilon<\min\left(\delta, (1-\delta^{\alpha_2})^{1/\alpha_1}\right)$
\BQN\label{new0}
 4\Pi_{\epsilon,1,\epsilon, \delta}(u)-\sum_{ 1\leq i<j\leq 4}P_{i,j}(u)\leq \Pi_{-1,1,-1,1}(u)&\leq& 4\left(\Pi_{0,1,0,\delta}(u)+\Pi_{0,(1-\delta^{\alpha_2})^{1/\alpha_1},0,1}(u)\right)
\EQN
where
$$P_{i,j}(u)=\pk{\sup_{(s,t)\in F^i}(B_{\alpha_1}(s)+B_{\alpha_2}(t))>u, \sup_{(s,t)\in F^j}(B_{\alpha_1}(s)+B_{\alpha_2}(t))>u},$$
with $F^1= F_{\epsilon,1,\epsilon, \delta}, \quad  F^2= F_{-1,-\epsilon,\epsilon, \delta}, \quad F^3= F_{\epsilon,1,-\delta, -\epsilon}, \quad F^4= F_{-1, -\epsilon, -\delta, -\epsilon}$.
Noting that for $i\neq j$,
\BQNY
\sup_{(s,t,s',t')\in F^i\times F^j}Var\left(B_{\alpha_1}(s)+B_{\alpha_2}(t)+B_{\alpha_1}(s')+B_{\alpha_2}(t')\right)<4-\delta_0,
\EQNY
by Borell-TIS theorem, we have
\BQNY
P_{i,j}(u)\leq \pk{\sup_{(s,t,s',t')\in F^i\times F^j}\left(B_{\alpha_1}(s)+B_{\alpha_2}(t)+B_{\alpha_1}(s')+B_{\alpha_2}(t')\right)>2u}
\leq e^{-\frac{(2u-a_0)^2}{2(4-\delta_0)}}=o\left(\Psi(u)\right), u\rw\IF,
\EQNY
with $a_0=2\E{\sup_{(s,t)E_{\alpha_1,\alpha_2}}B_{\alpha_1}+B_{\alpha_2}}<\IF$. Hence,
\BQN\label{new1}
\sum_{ 1\leq i<j\leq 4}P_{i,j}(u)=o\left(\Psi(u)\right), u\rw\IF.
\EQN

Let $g_1(t)=\frac{\alpha_1}{2}(1-|t|^{\alpha_2})^{1-1/\alpha_1}, \quad 0\leq t< 1$,  $g_2(t)=\frac{\alpha_2}{2}(1-|t|^{\alpha_1})^{1-1/\alpha_2}, \quad 0\leq t< 1$, $f_1(t)=(1-t^{\alpha_2})^{1/\alpha_1}, \quad 0\leq t< 1$,
 and $f_2(t)=(1-t^{\alpha_1})^{1/\alpha_2}, \quad 0\leq t< 1$.

\underline{Case $\alpha_2<1$}.  Since for  any $0<\delta_2<1$,
\BQN\label{new4}
0<\inf_{0\leq t\leq \delta_2}g_1(t)\leq \sup_{0\leq t\leq \delta_2}g_1(t)<\IF, \quad \inf_{0<t\leq \delta_2 }|f_1'(t)|=\inf_{0<t\leq \delta_2 }\left\{\frac{\alpha_2}{\alpha_1}(1-t^{\alpha_2})^{1/\alpha_1-1}t^{\alpha_2-1}\right\}>0,
\EQN
then in view of  (\ref{var6}) and (\ref{cor5}) and by Theorem \ref{th31} and Remark \ref{remark}, we have that,  for $0\leq \delta_1<\delta_2<1$,
\BQN\label{new2}
 \Pi_{0,1,\delta_1,\delta_2}(u)\sim \frac{2^{1-1/\alpha_1-1/\alpha_2}}{\alpha_1}\prod_{i=1}^2\mathcal{H}_{\alpha_i}\int_{\delta_1}^{\delta_2}(1-t^{\alpha_2})^{1/\alpha_1-1}dt
u^{2/\alpha_1+2/\alpha_2-2}\Psi(u).
\EQN
Exchanging the coordinates of $s$ and $t$ in (\ref{var5}) and (\ref{cor5}), we have that
\BQN\label{new}
1-\sigma(t,s)&\sim& \frac{\alpha_2}{2}(1-|t|^{\alpha_1})^{1-1/\alpha_2}\left||s|-(1-|t|^{\alpha_1})^{1/\alpha_2}\right|, \quad |t|^{\alpha_1}+|s|^{\alpha_2}\uparrow 1, |t|\leq 1-\delta,\nonumber\\
1-r(t,s,t_1,s_1))&\sim&\frac{|s-s_1|^{\alpha_2}}{2}+\frac{|t-t_1|^{\alpha_1}}{2}, \quad |s-s_1|, |t-t_1|\rw 0, |t|^{\alpha_1}+|s|^{\alpha_2}\uparrow 1.
\EQN
 Since for any $0<\delta_2<1$,
\BQN\label{new5}
0<\inf_{0\leq t\leq \delta_2}g_2(t)\leq \sup_{0\leq t\leq \delta_2}g_2(t)<\IF, \quad \inf_{0<t\leq \delta_2}|f_2'(t)|=\inf_{0<t\leq \delta_2}\left\{\frac{\alpha_1}{\alpha_2}(1-t^{\alpha_1})^{1/\alpha_2-1}t^{\alpha_1-1}\right\}>0,
\EQN
then in view of  (\ref{new}) and by Theorem \ref{th33} and Remark \ref{remark}, we have that,  for $0\leq \delta_1<\delta_2<1$,
\BQN\label{new3}
 \Pi_{\delta_1,\delta_2, 0,1}(u)\sim \frac{2^{1-1/\alpha_1-1/\alpha_2}}{\alpha_2}\prod_{i=1}^2\mathcal{H}_{\alpha_i}\int_{\delta_1}^{\delta_2}(1-t^{\alpha_1})^{1/\alpha_2-1}dt
u^{2/\alpha_1+2/\alpha_2-2}\Psi(u).
\EQN
Recalling that $0<\epsilon<\min\left(\delta, (1-\delta^{\alpha_2})^{1/\alpha_1}\right)$, we have
$$\Pi_{\epsilon,1,\epsilon, \delta}(u)\leq  \Pi_{0,1,\epsilon, \delta}(u)\leq \Pi_{\epsilon,1,\epsilon, \delta}(u)+\pk{\sup_{(s,t)\in E_{\alpha_1,\alpha_2}\setminus E(u)}(B_{\alpha_1}(s)+B_{\alpha_2}(t))>u}$$ with $E(u)=\{1- u^{-2}(\ln u)^2\leq |s|^{\alpha_1}+|t|^{\alpha_2}\leq 1\}$.
Using the fact that \BQNY
\sup_{(s,t)\in E_{\alpha_1,\alpha_2}\setminus E(u)}Var\left(B_{\alpha_1}(s)+B_{\alpha_2}(t)\right)&\leq& 1- u^{-2}(\ln u)^2,\\
Var\left(B_{\alpha_1}(s)+B_{\alpha_2}(t)-B_{\alpha_1}(s')-B_{\alpha_2}(t')\right)&\leq& 2\left(|s-s'|^{\alpha_1}+|t-t'|^{\alpha_1}\right), \quad (s,t), (s',t')\in E_{\alpha_1,\alpha_2},\EQNY
and by Theorem 8.1 in \cite{Pit96} (or Lemma 5.1 in \cite{KEP20161}), we have that, for $u$ sufficiently large,
\BQN\label{neg1}
\pk{\sup_{(s,t)\in E(u)}(B_{\alpha_1}(s)+B_{\alpha_2}(t))>u}\leq \mathbb{Q}u^{4/\alpha_1}\Psi\left(\frac{u}{\sqrt{1- u^{-2}(\ln u)^2}}\right)\leq \mathbb{Q}u^{4/\alpha_1}e^{-\frac{(\ln u)^2}{2}}\Psi\left(u\right),
\EQN
which combined with (\ref{new2}) leads to
\BQN\label{asym0}\Pi_{\epsilon,1,\epsilon, \delta}(u)\sim \Pi_{0,1,\epsilon, \delta}(u), \quad u\rw\IF.
\EQN
 Hence,
\BQNY
\limsup_{u\rw\IF}\frac{\Pi_{-1,1,-1,1}(u)}{2^{1-1/\alpha_1-1/\alpha_2}\prod_{i=1}^2\mathcal{H}_{\alpha_i}
u^{2/\alpha_1+2/\alpha_2-2}\Psi(u)}&\leq&  \frac{4}{\alpha_1}\int_{0}^{\delta}(1-t^{\alpha_2})^{1/\alpha_1-1}dt+ \frac{4}{\alpha_2}\int_{0}^{(1-\delta^{\alpha_2})^{1/\alpha_1}}(1-t^{\alpha_1})^{1/\alpha_2-1}dt,\\
\liminf_{u\rw\IF}\frac{\Pi_{-1,1,-1,1}(u)}{2^{1-1/\alpha_1-1/\alpha_2}\prod_{i=1}^2\mathcal{H}_{\alpha_i}
u^{2/\alpha_1+2/\alpha_2-2}\Psi(u)}&\geq&  \frac{4}{\alpha_1}\int_{\epsilon}^{\delta}(1-t^{\alpha_2})^{1/\alpha_1-1}dt.
\EQNY
Letting $\delta\rw 1$ in the above inequalities and by the fact that $0<\epsilon<\min\left(\delta, (1-\delta^{\alpha_2})^{1/\alpha_1}\right)$, we have
$$\Pi_{-1,1,-1,1}(u)\sim \frac{2^{3-1/\alpha_1-1/\alpha_2}}{\alpha_1}\prod_{i=1}^2\mathcal{H}_{\alpha_i} \int_{0}^{1}(1-t^{\alpha_2})^{1/\alpha_1-1}dt
u^{2/\alpha_1+2/\alpha_2-2}\Psi(u).$$
\underline{Case $\alpha_2=1$}. Note that if $\alpha_2=1$, (\ref{new4}) and (\ref{new5}) hold.
Thus in light of (\ref{var6}) and (\ref{cor5}) and by Theorem \ref{th31} and Remark \ref{remark}, we have that,  for $0\leq \delta_1<\delta_2<1$,
\BQNY
 \Pi_{0,1,\delta_1,\delta_2}(u)\sim \frac{2^{-1/\alpha_1}}{\alpha_1}\mathcal{H}_{\alpha_1}\int_{\delta_1}^{\delta_2}\widehat{\mathcal{P}}_1^{1}(1-t)^{1/\alpha_1-1}dt
u^{2/\alpha_1}\Psi(u).
\EQNY
By the fact that $\widehat{\mathcal{P}}_1^1=2$ (see, e.g., \cite{Pit96}), we have
\BQN\label{new6}
 \Pi_{0,1,\delta_1,\delta_2}(u)\sim  2^{1-1/\alpha_1}\left((1-\delta_1)^{1/\alpha_1}-(1-\delta_2)^{1/\alpha_1}\right)\mathcal{H}_{\alpha_1}
u^{2/\alpha_1}\Psi(u).
\EQN
In light of  (\ref{new}) and by Theorem \ref{th33} and Remark \ref{remark}, we have, for $0\leq \delta_1<\delta_2<1$,
\BQN\label{new7}
 \Pi_{\delta_1,\delta_2, 0,1}(u)\sim 2^{-1/\alpha_1}\mathcal{H}_{\alpha_1}\int_{\delta_1}^{\delta_2}\widehat{\mathcal{P}}_1^{1}dt
u^{2/\alpha_1}\Psi(u)\sim 2^{1-1/\alpha_1}(\delta_2-\delta_1)\mathcal{H}_{\alpha_1}
u^{2/\alpha_1}\Psi(u).
\EQN
Combining (\ref{new0})-(\ref{new1})and  (\ref{new6})-(\ref{new7}) and by the fact that $\Pi_{\epsilon,1,\epsilon, \delta}(u)\sim \Pi_{0,1,\epsilon, \delta}(u), u\rw\IF$( similarly as in (\ref{asym0})),  we have
\BQNY
\limsup_{u\rw\IF}\frac{\Pi_{-1,1,-1,1}(u)}{2^{1-1/\alpha_1}\mathcal{H}_{\alpha_1}
u^{2/\alpha_1}\Psi(u)}&\leq& 4\left(1- (1-\delta)^{1/\alpha_1}\right)+ 4(1-\delta^{\alpha_2})^{1/\alpha_1},\\
\liminf_{u\rw\IF}\frac{\Pi_{-1,1,-1,1}(u)}{2^{1-1/\alpha_1}\mathcal{H}_{\alpha_1}
u^{2/\alpha_1}\Psi(u)}&\geq&  4\left((1-\epsilon)^{1/\alpha_1}- (1-\delta)^{1/\alpha_1}\right).
\EQNY
Letting $\delta\rw 1$ in the above inequalities and noting that $0<\epsilon< (1-\delta)^{1/\alpha_1}$, we have
\BQNY
\Pi_{-1,1,-1,1}(u)\sim 2^{3-1/\alpha_1}\mathcal{H}_{\alpha_1}
u^{2/\alpha_1}\Psi(u).
\EQNY
\underline{Case $\alpha_2>1$}. In view of (\ref{new4}), we have that for any $0<\delta_1<\delta_2<1$,
\BQNY
0<\inf_{0\leq t\leq \delta_2}g_1(t)\leq \sup_{0\leq t\leq \delta_2}g_1(t)<\IF, \quad \inf_{\delta_1<t\leq \delta_2 }|f_1'(t)|=\inf_{\delta_1<t\leq \delta_2 }\left\{\frac{\alpha_2}{\alpha_1}(1-t^{\alpha_2})^{1/\alpha_2-1}t^{\alpha_2-1}\right\}>0.
\EQNY
Thus by Theorem \ref{th31} and Remark \ref{remark}, we have that,  for $0<\delta_1<\delta_2<1$,
\BQN\label{new8}
 \Pi_{0,1,\delta_1,\delta_2}(u)\sim \frac{ 2^{-1/\alpha_1}\alpha_2}{\alpha_1}\mathcal{H}_{\alpha_1}\int_{\delta_1}^{\delta_2}(1-t^{\alpha_2})^{1/\alpha_1-1}t^{\alpha_2-1}dt
u^{2/\alpha_1}\Psi(u).
\EQN
Note that if $\alpha_2>1$,  then $\lim_{t\rw 0}|f_1'(t)|=\lim_{t\rw 0}\frac{\alpha_2}{\alpha_1}(1-t^{\alpha_2})^{1/\alpha_2-1}t^{\alpha_2-1}=0$. This implies that this case isn't covered by Theorem \ref{th31}-\ref{th33}. We have to adopt another approach to deal with $ \Pi_{0,1,0,\delta_1}(u)$.
Observe that
$$ \Pi_{0,1,0,\delta_1}(u)\leq \pk{\sup_{(s,t)\in F(0,1, 0, \delta_1)\cap E(u)}(B_{\alpha_1}(s)+B_{\alpha_2}(t))>u}+ \pk{\sup_{(s,t)\in E_{\alpha_1,\alpha_2}\setminus E(u)}(B_{\alpha_1}(s)+B_{\alpha_2}(t))>u}.$$

Next we focus on the asymptotics of $\pk{\sup_{(s,t)\in F(0,1, 0, \delta_1)\cap E(u)}(B_{\alpha_1}(s)+B_{\alpha_2}(t))>u}$.
Let $I_{k,l}(u)=[ku^{-2/\alpha_1}, (k+1)u^{-2/\alpha_1}]\times[lu^{-2/\alpha_2}, (l+1)u^{-2/\alpha_2}]$. Observe that for $0<\delta_1<1/4$ and
$0\leq l\leq\left[\delta_1u^{2/\alpha_2}\right]+1$,
\BQNY
&&\{(s,t): lu^{-2/\alpha_2}\leq t\leq (l+1)u^{-2/\alpha_2}\}\cap\left(F(0,1, 0, \delta_1)\cap E(u)\right)\\
&&\quad \subset\left[\left(1-(l+1)^{\alpha_2}u^{-2}-u^{-2}(\ln u)^2\right)^{1/\alpha_1}, (1-l^{\alpha_2}u^{-2})^{1/\alpha_1}\right]\times [lu^{-2/\alpha_2}, (l+1)u^{-2/\alpha_2}]\\
&&\quad \subset\left[(1-l^{\alpha_2}u^{-2})^{1/\alpha_1}-b_l(u), (1-l^{\alpha_2}u^{-2})^{1/\alpha_1}\right]  \times [lu^{-2/\alpha_2}, (l+1)u^{-2/\alpha_2}],
\EQNY
with $b_l(u)=\frac{2}{\alpha_1}(1-\delta_1)^{(1/\alpha_1-1)\wedge 0}\left(\left((l+1)^{\alpha_2}-l^{\alpha_2}\right)u^{-2}+u^{-2}(\ln u)^2\right)$,  which implies that for any $0\leq l\leq\left[\delta_1u^{2/\alpha_2}\right]+1$ and $u$ large enough,
\BQNY
\#\{k: I_{k,l}(u)\cap \left(F(0,1, 0, \delta_1)\cap E(u)\right)\neq \emptyset\}\leq u^{2/\alpha_1}b_l(u)+1
\EQNY
Hence, by (\ref{cor5}) and Lemma \ref{PIPI},
\BQNY
&&\pk{\sup_{(s,t)\in F(0,1, 0, \delta_1)\cap E(u)}(B_{\alpha_1}(s)+B_{\alpha_2}(t))>u}\\
&&\quad \leq \sum_{l=0}^{\left[\delta_1u^{2/\alpha_2}\right]+1}\sum_{I_{k,l}(u)\cap F(0,1, 0, \delta_1)\cap E(u)\neq \emptyset}\pk{\sup_{(s,t)\in I_{k,l}(u)}\overline{B_{\alpha_1}(s)+B_{\alpha_2}(t)}>u}\\
&&\quad \leq \sum_{l=0}^{\left[\delta_1u^{2/\alpha_2}\right]+1}\left(u^{2/\alpha_1}b_l(u)+1\right)\prod_{i=1}^2\mathcal{H}_{\alpha_i}[0,2^{-1/\alpha_i}]\Psi(u)\\
&& \quad \leq  \prod_{i=1}^2\mathcal{H}_{\alpha_i}[0,2^{-1/\alpha_i}]\left(\frac{4}{\alpha_1}(1-\delta_1)^{(1/\alpha_1-1)\wedge 0}\left(\delta_1^{\alpha_2}u^{2/\alpha_1}+\delta_1u^{2/\alpha_1+2/\alpha_2-2}(\ln u)^2\right)+2\delta_1u^{2/\alpha_2}\right)\Psi(u)\\
&&\quad \leq \mathbb{Q}\delta_1u^{2/\alpha_1}\Psi(u), \quad u\rw\IF,
\EQNY
which combined with (\ref{neg1}) yields that
\BQN\label{new9}
\limsup_{u\rw\IF}\frac{\Pi_{0,1,0,\delta_1}(u)}{u^{2/\alpha_1}\Psi(u)}\leq \mathbb{Q}\delta_1.
\EQN
Next we focus on $\Pi_{0,\delta_2,0,1}(u)$. Observe that, for any $0<\delta_2<1/4$,
$$ \Pi_{0,\delta_2,0,1}(u)\leq \pk{\sup_{(s,t)\in F(0,\delta_2, 0, 1)\cap  E(u)}(B_{\alpha_1}(s)+B_{\alpha_2}(t))>u}+ \pk{\sup_{(s,t)\in E_{\alpha_1,\alpha_2}\setminus E(u)}(B_{\alpha_1}(s)+B_{\alpha_2}(t))>u},$$
 and for
$0\leq k\leq\left[\delta_2u^{2/\alpha_1}\right]+1$,
\BQNY
&&\{(s,t): ku^{-2/\alpha_1}\leq s\leq (k+1)u^{-2/\alpha_1}\}\cap\left(F(0,1, 0, \delta_1)\cap E(u)\right)\\
&&\quad \subset[ku^{-2/\alpha_1}, (k+1)u^{-2/\alpha_1}]\times\left[\left(1-(k+1)^{\alpha_1}u^{-2}-u^{-2}(\ln u)^2\right)^{1/\alpha_2}, (1-k^{\alpha_1}u^{-2})^{1/\alpha_2}\right]\\
&&\quad \subset[ku^{-2/\alpha_1}, (k+1)u^{-2/\alpha_1}]\times\left[(1-k^{\alpha_1}u^{-2})^{1/\alpha_2}-\widetilde{b}_{k}(u), (1-k^{\alpha_1}u^{-2})^{1/\alpha_2}\right],
\EQNY
with $\widetilde{b}_{k}(u)=\frac{2}{\alpha_2}(1-\delta_2)^{1/\alpha_2-1}\left(\left((k+1)^{\alpha_1}-k^{\alpha_1}\right)u^{-2}+u^{-2}(\ln u)^2\right)$.
Since,  for
$0\leq k\leq\left[\delta_2u^{2/\alpha_1}\right]+1$, with $\theta\in (0,1)$,
\BQNY
u^{2/\alpha_2}\widetilde{b}_{k}(u)&=&\frac{2}{\alpha_2}(1-\delta_2)^{1/\alpha_2-1}\left(\alpha_1(k+\theta)^{(\alpha_1-1)\vee 0}u^{2/\alpha_2-2}+u^{2/\alpha_2-2}(\ln u)^2\right)\\
&\leq& \mathbb{Q}\left(u^{\min(2/\alpha_2-2/\alpha_1, 2/\alpha_2-2)} +u^{2/\alpha_2-2}(\ln u)^2\right)\rw 0, \quad u\rw\IF,
\EQNY
then for $0\leq k\leq\left[\delta_2u^{2/\alpha_1}\right]+1$ and $u$ sufficiently large,
$$\#\{l: I_{k,l}(u)\cap F(0,1, 0, \delta_1)\cap E(u)\neq \emptyset\}\leq 2.$$
Thus by (\ref{cor5}) and Lemma \ref{PIPI}, we have
\BQNY
&&\pk{\sup_{(s,t)\in F(0,\delta_2, 0, 1)\cap E(u)}(B_{\alpha_1}(s)+B_{\alpha_2}(t))>u}\\
&&\quad \leq \sum_{k=0}^{\left[\delta_2u^{2/\alpha_1}\right]+1}\sum_{I_{k,l}(u)\cap F(0,\delta_2, 0,1)\cap E(u)\neq \emptyset}\pk{\sup_{(s,t)\in I_{k,l}(u)}\overline{B_{\alpha_1}(s)+B_{\alpha_2}(t)}>u}\\
&&\quad \leq 2\sum_{k=0}^{\left[\delta_2u^{2/\alpha_1}\right]+1}\prod_{i=1}^2\mathcal{H}_{\alpha_i}[0,2^{-1/\alpha_i}]\Psi(u)\\
&& \quad \leq 4\delta_2\prod_{i=1}^2\mathcal{H}_{\alpha_i}[0,2^{-1/\alpha_i}]u^{2/\alpha_1}\Psi(u),
\EQNY
which together with (\ref{neg1}) implies that
\BQN\label{new10}
\limsup_{u\rw\IF}\frac{\Pi_{0,\delta_2,0,1}(u)}{u^{2/\alpha_1}\Psi(u)}\leq 4\delta_2\prod_{i=1}^2\mathcal{H}_{\alpha_i}[0,2^{-1/\alpha_i}]\leq \mathbb{Q}\delta_2.
\EQN
By the fact that $\Pi_{\epsilon,1,\epsilon, \delta}(u)\sim \Pi_{0,1,\epsilon, \delta}(u), u\rw\IF$ (similarly as given in (\ref{asym0})), and $\Pi_{0,1,0, \delta}(u)\leq \Pi_{0,1,\epsilon, \delta}(u)+\Pi_{0,1,0,\epsilon}(u)$,  combination of  (\ref{new0})-(\ref{new1}) and (\ref{new8})-(\ref{new10}) leads to
\BQNY
\limsup_{u\rw\IF}\frac{\Pi_{-1,1,-1,1}(u)}{u^{2/\alpha_1}\Psi(u)}&\leq& \frac{ 2^{2-1/\alpha_1}\alpha_2}{\alpha_1}\mathcal{H}_{\alpha_1}\int_{\epsilon}^{\delta}(1-t^{\alpha_2})^{1/\alpha_1-1}t^{\alpha_2-1}dt+
\mathbb{Q}(\epsilon+(1-\delta^{\alpha_2})^{1/\alpha_1}),\\
\liminf_{u\rw\IF}\frac{\Pi_{-1,1,-1,1}(u)}{u^{2/\alpha_1}\Psi(u)}&\geq& \frac{ 2^{2-1/\alpha_1}\alpha_2}{\alpha_1}\mathcal{H}_{\alpha_1}\int_{\epsilon}^{\delta}(1-t^{\alpha_2})^{1/\alpha_1-1}t^{\alpha_2-1}dt.
\EQNY
Letting $\delta\rw 1$ in the above inequalities and recalling that $0<\epsilon<\min(\delta, (1-\delta^{\alpha_2})^{1/\alpha_1})$, we have
\BQNY
\Pi_{-1,1,-1,1}(u)\sim \frac{ 2^{2-1/\alpha_1}\alpha_2}{\alpha_1}\mathcal{H}_{\alpha_1}\int_{0}^{1}(1-t^{\alpha_2})^{1/\alpha_1-1}t^{\alpha_2-1}dt u^{2/\alpha_1}\Psi(u).
\EQNY
This completes the proof. \QED

\proofprop{cor12} Let $\alpha_1=\alpha_2=\alpha$.   Note that in this case, (\ref{var5})-(\ref{new1}) hold with $\alpha_1, \alpha_2$ replaced by $\alpha$. We use the same notation as  in the proof of Corollary \ref{cor11}.\\
\underline{Case $\alpha<1$}. Using the same arguments as in the proof of Corollary \ref{cor11} and by Theorem \ref{th32} and Remark \ref{remark}, we have that
$$\Pi_{-1,1,-1,1}(u)\sim \frac{2^{3-2/\alpha}}{\alpha}\left(\mathcal{H}_{\alpha}\right)^2 \int_{0}^{1}(1-t^{\alpha})^{1/\alpha-1}dt
u^{4/\alpha-2}\Psi(u).$$
\underline{Case $\alpha=1$}. Note that $1-\sigma(s,t)
 \sim \frac{1}{2}|s-(1-t)|$ as $s+t\uparrow 1$ for $0\leq s\leq 1, t\geq 0$. $f_1(t)=1-t$ and $g_1(t)=1/2$. $|f'(t)|=1, 0<t<1$.
 In light of  (\ref{cor5}) and by Theorem \ref{th32} and Remark \ref{remark}, we have that,  for $0\leq \delta_1<\delta_2\leq 1$,
\BQNY
 \Pi_{0,1,\delta_1,\delta_2}(u)\sim 2^{-1}\int_{\delta_1}^{\delta_2}\widehat{\mathcal{H}}_1^{1,-1}dt
u^{2}\Psi(u)\sim 2^{-1}(\delta_2-\delta_1)\widehat{\mathcal{H}}_{1}^{1,-1}u^2\Psi(u).
\EQNY
Similarly,
$$
\Pi_{\delta_1,\delta_2, 0,1}(u)\sim 2^{-1}\int_{\delta_1}^{\delta_2}\widehat{\mathcal{H}}_1^{1,-1}dt
u^{2}\Psi(u)\sim 2^{-1}(\delta_2-\delta_1)\widehat{\mathcal{H}}_{1}^{1,-1}u^2\Psi(u).$$
Following same arguments as  in the proof of Corollary \ref{cor11}, we have that
$$\Pi_{-1,1,-1,1}(u)\sim 2\widehat{\mathcal{H}}_{1}^{1,-1}u^2\Psi(u).$$
\underline{Case $\alpha>1$}. Using Theorem \ref{th32} and Remark \ref{remark}, we have that,  for $0<\delta_1<\delta_2<1$,
\BQNY
 \Pi_{0,1,\delta_1,\delta_2}(u)\sim 2^{-1/\alpha}\mathcal{H}_{\alpha}\int_{\delta_1}^{\delta_2}\left((1-t^{\alpha})^{1-\alpha}t^{\alpha(\alpha-1)}+1\right)^{1/\alpha}dt
u^{2/\alpha}\Psi(u).
\EQNY
Following same arguments as in (\ref{new9})-(\ref{new10}), we have that
$$\lim_{\delta_1\rw 0}\lim_{u\rw\IF}\frac{ \Pi_{0,1,0,\delta_1}(u)}{u^{2/\alpha}\Psi(u)}=\lim_{\delta_2\rw 0}\lim_{u\rw\IF}\frac{ \Pi_{0,\delta_2,0,1}(u)}{u^{2/\alpha}\Psi(u)}=0.$$
Hence,
\BQNY
 \Pi_{0,1,0,1}(u)\sim 2^{-1/\alpha}\mathcal{H}_{\alpha}\int_{0}^{1}\left((1-t^{\alpha})^{1-\alpha}t^{\alpha(\alpha-1)}+1\right)^{1/\alpha}dt
u^{2/\alpha}\Psi(u).
\EQNY
This completes the proof. \QED
 \section{Appendix}
 Following Theorem 2.1 in \cite{KEP20161} and Lemma 5.2 in \cite{KEP20162}, we present the uniform expansion of tail asymptotics of supremum of Gaussian fields over short intervals. \\
 Let $\rho_i\in \mathcal{R}_{\alpha_i/2}$, $v_i\in \mathcal{R}_{\beta_i/2}, i=1,2$ be non-negative functions with $0<\alpha_i\leq 2, \beta_i>0, i=1,2.$
Let $X_{u,k}(s,t), k\in \mathcal{K}_u$, with $\mathcal{K}_u$ representing the index set, be  centered Gaussian random fields over $\mathcal{E}(u):=\{((1+o(1))\overleftarrow{\LL_1}(u^{-1})s, (1+o(1))\overleftarrow{\LL_2}(u^{-1}) t ) , (s,t)\in \mathcal{E}\}$ with $\mathcal{E}$ an compact set containing $0$. Moreover, $h_k(u), k\in K_u$ are positive functions of $u$ satisfying $\lim_{u\rw\IF}\frac{h_k(u)}{u}=1$ uniformly with respect to $k\in K_u$.
Suppose further that $X_{u,k}$ has unit variance, continuous trajectories  and correlation function $r_k(s,t,s_1,t_1)$ satisfying (\ref{cor2}) uniformly with respect to $(k,l)\in K_u$.
\BEL\label{PIPI}
Let  $d_u(s,t),u>0$ be continuous functions satisfying
 \BQN\label{eq1}
 \lim_{u\rw\IF}\sup_{(s,t)\in\mathcal{E}(u), k\in \mathcal{K}_u}\left|h_k^2(u)d_u(\overleftarrow{\LL_1}(u^{-1})s, \overleftarrow{\LL_2}(u^{-1})t)-d(s,t)\right|=0.
 \EQN
  Then 
\BQNY
\lim_{u\rw\IF}\sup_{k\in \mathcal{K}_u}\left|(\Psi(h_k(u)))^{-1}\mathbb{P}\left(\sup_{(s,t)\in \mathcal{E}(u)}\frac{X_{u,k}(s,t)}{1+d_u(s,t)}>h_k(u)\right)-
\E{e^{\sup_{(s,t)\in \mathcal{E}}\{W_{\alpha_1,\alpha_2}(s,t)-d(s,t)\}}}\right|=0,
\EQNY
with $W_{\alpha_1, \alpha_2}$ defined right before (\ref{piterbarg}).
\EEL
{\bf Acknowledgement}: Thanks to Krzysztof D\c{e}bicki and  Enkelejd Hashorva for their suggestions.  Thanks to Swiss National Science Foundation grant No. 200021-166274.
\bibliographystyle{ieeetr}

 \bibliography{sigmaABCDEDD}

\end{document}